\documentclass{article}
 \usepackage{epsfig}\usepackage{graphicx}
\usepackage[latin1]{inputenc}
\usepackage{amsmath}\usepackage{amsthm}\usepackage{amsfonts}
\usepackage{amssymb}\usepackage{amstext}\usepackage{amsgen}
\usepackage{amsbsy}\usepackage{amsopn}
\usepackage{hyperref}
\usepackage{mathrsfs, a4wide, color}
\usepackage{pifont}

\numberwithin{equation}{section}

\theoremstyle{plain}
\newtheorem{theorem}{Theorem}[section]
\newtheorem{corollary}[theorem]{Corollary}

\newtheorem{lemma}[theorem]{Lemma}

\theoremstyle{definition}
\newtheorem{definition}[theorem]{Definition}

\newtheorem{remark}[theorem]{Remark}

\newcommand{\e}{\ensuremath{\mathrm{e\;\!}}}
\newcommand{\im}{\ensuremath{\mathrm{i}}}

\newcommand{\Var}{\mathop{\mathrm{Var_{\mathbb
Q^S}}}}
\newcommand{\Ex}{{\mathbb E}}
\newcommand{\Pa}{{\mathbb P}}
\newcommand{\Q}{{\mathbb Q}}
\newcommand{\C}{{\mathbb C}}
\newcommand{\R}{{\mathbb R}}
\newcommand{\N}{{\mathbb N}}

\newcommand{\Bcal}{{\mathcal B}}

\newcommand{\Dcal}{{\mathcal D}}
\newcommand{\Ecal}{{\mathcal E}}
\newcommand{\Fcal}{{\mathcal F}}

\newcommand{\Scal}{{\mathcal S}}

\newcommand{\Wcal}{{\mathcal W}}
\newcommand{\Xcal}{{\mathcal X}}

\begin{document}
\title{Affine Diffusion Processes: Theory and Applications}
\author{Damir Filipovi\'c \& Eberhard Mayerhofer
\footnote{University of Vienna, and Vienna University of Economics
and Business Administration, Heiligenst\"adter Stra{\ss}e 46-48,
1190 Vienna, Austria, email: \{damir.filipovic,
eberhard.mayerhofer\}@vif.ac.at. We thank Paul Glasserman for
helpful comments. Financial support from WWTF (Vienna Science and
Technology Fund) is gratefully acknowledged.} }

\date{20 February, 2009}
\maketitle \thispagestyle{empty}

\abstract{We revisit affine diffusion processes on general and on
the canonical state space in particular. A detailed study of
theoretic and applied aspects of this class of Markov processes is
given. In particular, we derive admissibility conditions and provide
a full proof of existence and uniqueness through stochastic
invariance of the canonical state space. Existence of exponential
moments and the full range of validity of the affine transform
formula are established. This is applied to the pricing of bond and
stock options, which is illustrated for the Vasi\v{c}ek,
Cox--Ingersoll--Ross and Heston models.}

\section{Introduction} Affine Markov models have been employed in
finance since decades, and they have found growing interest due to
their computational tractability as well as their capability to
capture empirical evidence from financial time series. Their main
applications lie in the theory of term structure of interest rates,
stochastic volatility option pricing and the modeling of credit risk
(see \cite{duffilsch:03} and the references therein). There is a
vast literature on affine models. We mention here explicitly just
the few articles \cite{andersenpiterbarg, brownschaefer,
chefilkim08, daisin00, dufkan, dufpansing, Glasserman, hes93,
kellermoments, keller, Lee} and \cite{duffilsch:03} for a broader
overview.

In this paper, we revisit the class of affine diffusion processes on
subsets of $\mathbb R^d$ and on the canonical state space $\mathbb
R_+^m\times \mathbb R^n$, in particular. In
Section~\ref{secdefchaff}, we first provide necessary and sufficient
conditions on the parameters of a diffusion process $X$ to satisfy
the affine transform formula
\[\Ex\left[ \e^{u^\top  X(T)}\mid\Fcal_t\right] =
\e^{\phi(T-t,\,u)+\psi(T-t,\,u)^\top   X(t)}.\] The functions $\phi$
and $\psi$ in turn are given as solutions of a system of coupled
Riccati equations. Arguing by stochastic invariance, in
Section~\ref{seccanstatespace}, we can further restrict the choice
of admissible diffusion parameters.

Glasserman and Kim~\cite{Glasserman} showed recently that the affine
transform formula holds whenever either side is well defined under
the assumption of strict mean reversion. This is an extension of the
findings in \cite{duffilsch:03}, where only sufficient conditions
are given in terms of analyticity of the right hand side. The strict
mean reversion assumption, however, excludes the Heston stochastic
volatility model. In our paper, we show that strict mean reversion
is not needed (Theorem~\ref{thmextanaXXX}). As a by product, we
obtain some non-trivial convexity results for Riccati equations.
Having the full range of validity of the above transform formula
under control, in Section~\ref{secdiscpri}, we can then proceed to
pricing bond and stock options in affine models. Particular examples
are the Vasi\v{c}ek and Cox--Ingersoll--Ross (CIR) short rate models
in Section~\ref{secbondopaffine}, and Heston's stochastic volatility
model in Section~\ref{secheston}.

The representation of affine short rate models bears some ambiguity
with respect to linear transformations of the state process. This
motivates the question whether there exists a classification method
ensuring that affine short rate models with the same observable
implications have a unique canonical representation. This topic has
been addressed in \cite{daisin00,coldufetal06,jos06,chefilkim08}. In
Section~\ref{secafftransf}, we recap this issue and show that the
diffusion matrix of $X$ can always be brought into block-diagonal
form by a regular linear transform leaving the canonical state space
invariant.

The existence and uniqueness question of the relevant stochastic
differential equation is completely solved through stochastic
invariance and the block-diagonal transformation in
Section~\ref{secexiuniaffine}. The presented proof builds on the
seminal result by Yamada and Watanabe~\cite{yamwat71}. We therefore
approach the existence issue differently from \cite{duffilsch:03}
which uses infinite divisibility on the canonical state space and
the Markov semigroup theory.

In the appendix, we provide some self contained proofs of existence
and comparison statements for relevant systems of Riccati equations
(Section~\ref{sec: comp}). Moreover, some moment lemmas from
\cite{duffilsch:03} in a more elaborated fashion can be found in
Section~\ref{secregcharfun}.

\section{Definition and Characterization of Affine
Processes}\label{secdefchaff} Fix a dimension $d\ge 1$ and a closed
state space $\Xcal\subset\mathbb R^d$ with non-empty interior. We
let $b:\Xcal\to\mathbb R^d$ be continuous, and $\rho:\Xcal\to\mathbb
R^{d\times d}$ be measurable and such that the diffusion matrix
\[ a(x)=\rho(x)\rho(x)^\top \]
is continuous in $x\in\Xcal$. Let $W$ denote a $d$-dimensional
Brownian motion defined on a filtered probability space
$(\Omega,\Fcal,(\Fcal_t),\Pa)$. Throughout, we assume that for every
$x\in\Xcal$ there exists a unique solution $X=X^x$ of the stochastic
differential equation
\begin{equation}\label{sdeaffgen}
dX(t)=b(X(t))\,dt+\rho(X(t))\,dW(t),\quad X(0)={x}.
\end{equation}
\begin{definition}
We call $X$ {\emph{affine}} if the $\Fcal_t$-conditional
characteristic function of $X(T)$ is exponential affine in $X(t)$,
for all $t\le T$. That is, there exist $\mathbb C$- and $\mathbb
C^d$-valued functions $\phi(t,u)$ and $\psi(t,u)$, respectively,
with jointly continuous $t$-derivatives such that $X=X^x$ satisfies
\begin{equation}\label{eqaffdef1}
 \Ex\left[ \e^{u^\top  X(T)}\mid\Fcal_t\right] =
\e^{\phi(T-t,\,u)+\psi(T-t,\,u)^\top   X(t)}
\end{equation}
for all $u\in\im\mathbb R^d$, $t\le T$ and $x\in\Xcal$.
\end{definition}

Since the conditional characteristic function is bounded by one, the
real part of the exponent $\phi(T-t,u)+\psi(T-t,u)^\top   X(t)$ in
\eqref{eqaffdef1} has to be negative. Note that $\phi(t,u)$ and
$\psi(t,u)$ for $t\ge 0$ and $u\in\im\mathbb R^d$ are
uniquely\footnote{In fact, $\phi(t,u)$ may be altered by multiples
of $2\pi\im$. We uniquely fix the continuous function $\phi(t,u)$ by
the initial condition $\phi(0,u)=0$.} determined by
\eqref{eqaffdef1}, and satisfy the initial conditions $\phi(0,u)=0$
and $\psi(0,u)=u$, in particular.

We first derive necessary and sufficient conditions for $X$ to be
affine.
\begin{theorem}\label{thmaffchar}
Suppose $X$ is affine. Then the diffusion matrix $a(x)$ and drift
$b(x)$ are affine in $x$. That is,
  \begin{equation}\label{eqaffba}
  \begin{aligned}
    a(x)&=a+\sum_{i=1}^d x_i \alpha_i\\
    b(x)&=b+\sum_{i=1}^d x_i \beta_i=b+\Bcal x\\
  \end{aligned}
  \end{equation}
for some $d\times d$-matrices $a$ and $\alpha_i$, and $d$-vectors
$b$ and $\beta_i$, where we denote by
\[\Bcal=(\beta_1,\dots,\beta_d)\] the $d\times d$-matrix with
$i$-th column vector $\beta_i$, $1\le i\le d$. Moreover, $\phi$ and
$\psi=(\psi_1,\dots,\psi_d)^\top$ solve the system of Riccati
equations
  \begin{equation}\label{eqriccfull}
    \begin{aligned}
    \partial_t\phi(t,u)&=\frac{1}{2}\psi(t,u)^\top
   a \,\psi(t,u)+b^\top \psi(t,u)\\
   \phi(0,u)&=0\\
      \partial_t\psi_i(t,u)&=\frac{1}{2}\psi(t,u)^\top
    \alpha_i\,\psi(t,u)+\beta_i^\top
      \psi(t,u) ,\quad 1\le i\le d,\\
    \psi(0,u)&=u.
    \end{aligned}
  \end{equation}
In particular, $\phi$ is determined by $\psi$ via simple
integration:
\[ \phi(t,u)=\int_0^t \left(\frac{1}{2}\psi(s,u)^\top
   a \,\psi(s,u)+b^\top \psi(s,u)\right)ds.\]

Conversely, suppose the diffusion matrix $a(x)$ and drift $b(x)$ are
affine of the form \eqref{eqaffba} and suppose there exists a
solution $(\phi,\psi)$ of the Riccati equations~\eqref{eqriccfull}
such that $\phi(t,u)+\psi(t,u)^\top x$ has negative real part for
all $t\ge 0$, $u\in\im\mathbb R^d$ and $x\in\Xcal$. Then $X$ is
affine with conditional characteristic function~\eqref{eqaffdef1}.
\end{theorem}

\begin{proof}

Suppose $X$ is affine. For $T>0$ and $u\in\im\mathbb R^d$ define the
complex-valued It\^o process
\[M(t)=\e^{\phi(T-t,u)+\psi(T-t,u)^\top X(t)}.\] We can apply It\^o's
formula, separately to real and imaginary part of $M$, and obtain
\[ dM(t)=I(t)\,dt +\psi(T-t,u)^\top\rho(X(t))\,dW(t),\quad t\le T,\]
with
\begin{multline*}
I(t)= -\partial_T\phi(T-t,u)-\partial_T\psi(T-t,u)^\top
X(t)\\+\psi(T-t,u)^\top b(X(t))+\frac{1}{2}\psi(T-t,u)^\top
a(X(t))\,\psi(T-t,u).
\end{multline*}
Since $M$ is a martingale, we have $I(t)=0$ for all $t\le T$ a.s.
Letting $t\to 0$, by continuity of the parameters, we thus obtain
\[\partial_T\phi(T,u)+\partial_T\psi(T,u)^\top x=\psi(T,u)^\top
b(x)+\frac{1}{2}\psi(T,u)^\top a(x)\,\psi(T,u)\] for all
$x\in\Xcal$, $T\ge 0$, $u\in\im\mathbb R^d$. Since $\psi(0,u)=u$,
this implies that $a$ and $b$ are affine of the form
\eqref{eqaffba}. Plugging this back into the above equation and
separating first order terms in $x$ yields \eqref{eqriccfull}.

Conversely, suppose $a$ and $b$ are of the form \eqref{eqaffba}. Let
$(\phi,\psi)$ be a solution of the Riccati equations
\eqref{eqriccfull} such that $\phi(t,u)+\psi(t,u)^\top x$ has
negative real part for all $t\ge 0$, $u\in\im\mathbb R^d$ and
$x\in\Xcal$. Then $M$, defined as above, is a uniformly bounded
local martingale, and hence a martingale, with $M(T)=\e^{u^\top
X(T)}$. Therefore $\Ex[M(T)\mid\Fcal_t]=M(t)$, for all $t\le T$,
which is \eqref{eqaffdef1}, and the theorem is proved.
\end{proof}

We now recall an important global existence, uniqueness and
regularity result for the above Riccati equations. We let $K$ be a
placeholder for either $\mathbb R$ or $\mathbb C$.

\begin{lemma}\label{lemextanalytic1}
Let $a$ and $\alpha_i$ be real $d\times d$-matrices, and $b$ and
$\beta_i$ be real $d$-vectors, $1\le i\le d$.

\begin{enumerate}
  \item\label{lemextanalytic11} For every $u\in K^d$, there exists
  some $t_+(u)\in (0,\infty]$ such that there exists a unique solution
$(\phi(\cdot,u),\psi(\cdot,u)):[0,t_+(u))\to K\times K^d$ of the
Riccati equations \eqref{eqriccfull}. In particular,
$t_+(0)=\infty$.

\item\label{lemextanalytic12} The domain \[\Dcal_K=\{(t,u)\in\R_+\times K^d\mid t<t_+(u)\}\]
is open in $\mathbb R_+\times K^d$ and maximal in the sense that for
all $u\in K^d$ either $t_+(u)=\infty$ or $\lim_{t\uparrow t_+(u)} \|
\psi(t,u)\|=\infty$, respectively, .

\item\label{lemextanalytic13} For every $t\in\mathbb R_+$, the $t$-section
\[\Dcal_K(t)=\{ u\in K^d\mid (t,u)\in\Dcal_K\}\] is an open
neighborhood of $0$ in $K^d$. Moreover, $\Dcal_K(0)=K^d$ and
$\Dcal_K(t_1)\supseteq\Dcal_K(t_2)$ for $0\le t_1\le t_2$.

\item\label{lemextanalytic14} $\phi$ and $\psi$ are analytic functions on $\Dcal_K$.

\item\label{lemextanalytic15} $\Dcal_{\mathbb R}=\Dcal_\mathbb C\cap (\R_+\times\R^d)$.
\end{enumerate}
\end{lemma}

Henceforth, we shall call $\Dcal_K$ the maximal domain for equation
\eqref{eqriccfull}.

\begin{proof}
Since the right-hand side of \eqref{eqriccfull} is formed by
analytic functions in $\psi$ on $K^d$, part~\ref{lemextanalytic11}
follows from the basic theorems for ordinary differential equations,
e.g.\ \cite[Theorem 7.4]{ama90}. In particular, $t_+(0)=\infty$
since $(\phi(\cdot,0),\psi(\cdot,0))\equiv 0$ is the unique solution
of \eqref{eqriccfull} for $u=0$. It is proved in \cite[Theorems 7.6
and 8.3]{ama90} that $\Dcal_K$ is maximal and open, which is
part~\ref{lemextanalytic12}. This also implies that all $t$-sections
$\Dcal_K(t)$ are open in $K^d$. The inclusion
$\Dcal_K(t_1)\supseteq\Dcal_K(t_2)$ is a consequence of the
maximality property from part~\ref{lemextanalytic12}. Whence
part~\ref{lemextanalytic13} follows. For a proof of
part~\ref{lemextanalytic14} see \cite[Theorem 10.8.2]{dieu60}.
Part~\ref{lemextanalytic15} is obvious.
\end{proof}

We will provide in Section~\ref{sec: comp} below some substantial
improvements of the properties stated in Lemma~\ref{lemextanalytic1}
for the canonical state space $\Xcal$ introduced in the following
section.

\section{Canonical State Space}\label{seccanstatespace}

There is an implicit trade off between the parameters
$a,\,\alpha_i,\,b,\,\beta_i$ in \eqref{eqaffba} and the state space
$\Xcal$:
\begin{itemize}
  \item $a,\,\alpha_i,\,b,\,\beta_i$ must be such that $X$ does not
  leave the set $\Xcal$, and
  \item $a,\,\alpha_i$ must be such that $a+\sum_{i=1}^d
  x_i\alpha_i$ is symmetric and positive semi-definite for all
  $x\in\Xcal$.
\end{itemize}
To gain further explicit insight into this interplay, we now and
henceforth assume that the state space is of the following canonical
form
 \[\Xcal=\mathbb R^m_+\times\mathbb R^n\]
 for some integers $m,n\ge 0$ with $m+n=d$.

\begin{remark}\label{rem: degenerate}
This canonical state space covers most applications appearing in the
finance literature. However, other choices for the state space of an
affine process are possible:

\begin{enumerate}
\item For instance, the following example for $d=1$ admits as state
space any closed interval $\Xcal\subset\mathbb R$ containing 0:
\[ dX=-X\,dt,\quad X(0)=x\in\Xcal.\]
This degenerate diffusion process is affine, since $\e^{u
X(T)}=\e^{u\e^{-(T-t)}X(t)} $ for all $t\le T$ (\cite{duffilsch:03},
Section 12). In general, affine diffusion processes on compact state
spaces have to be degenerate.
\item Matrix state-spaces $S_d^+$ ($d\geq 2$), the cone of symmetric
positive definite matrices (see \cite{Bru,BPT,TFG,CFMT,TG}.
\item Parabolic state-spaces, cf.\ \cite{sufana}, which are in turn, related
to quadratic processes on the canonical state-space
(\cite{chenfilipoor}, see also their Example 5.3, section 5)
\end{enumerate}
\end{remark}

For the above canonical state space, we can give necessary and
sufficient admissibility conditions on the parameters. The following
terminology will be useful in the sequel. We define the index sets
\[{I}=\{1,\dots,m\}\quad\text{and}\quad {J}=\{m+1,\dots,m+n\}.\]
For any vector $\mu$ and matrix $\nu$, and index sets $M,N$, we
denote by
\[ \mu_M=(\mu_i)_{i\in M},\quad
\nu_{MN}=(\nu_{ij})_{i\in M,\,j\in N}\] the respective sub-vector
and -matrix.

\begin{theorem}\label{thmaffine}
The process $X$ on the canonical state space $\mathbb
R^m_+\times\mathbb R^n$ is affine if and only if $a(x)$ and $b(x)$
are affine of the form \eqref{eqaffba} for parameters
$a,\,\alpha_i,\,b,\,\beta_i$ which are admissible in the following
sense:
\begin{equation}\label{eqadmiss}
\begin{aligned}
a,\,\alpha_i &\;\text{ are symmetric positive semi-definite,}\\
a_{II}&=0\quad\text{(and thus $a_{IJ}=a_{JI}^\top=0$),}\\ \alpha_j
&=0\quad\text{for all}\quad j\in J\\
\alpha_{i,kl}=\alpha_{i,lk}&=0\quad\text{for $k\in I\setminus\{i\}$,
for all $1\le i,l\le d$,}\\ b&\in \mathbb R^m_+\times \mathbb R^n,\\
\Bcal_{IJ}&=0,\\ \Bcal_{II}&\;\text{ has positive off-diagonal
elements.}
\end{aligned}
\end{equation}
In this case, the corresponding system of Riccati
equations~\eqref{eqriccfull} simplifies to
  \begin{equation}\label{eqriccfulladmiss}
    \begin{aligned}
    \partial_t\phi(t,u)&=\frac{1}{2}\psi_J(t,u)^\top
   a_{JJ} \,\psi_J(t,u)+b^\top \psi(t,u)\\
   \phi(0,u)&=0\\
      \partial_t\psi_i(t,u)&=\frac{1}{2}\psi(t,u)^\top
    \alpha_i\,\psi(t,u)+\beta_i^\top
      \psi(t,u),\quad i\in I,\\
\partial_t\psi_J(t,u)&=\Bcal^\top_{JJ}
      \psi_J(t,u),\\
    \psi(0,u)&=u,
    \end{aligned}
  \end{equation}
and there exists a unique global solution
$(\phi(\cdot,u),\psi(\cdot,u)):\mathbb R_+\to\mathbb
C_-\times\mathbb C^m_-\times\im\mathbb R^n$ for all initial values
$u\in\mathbb C^m_-\times\im\mathbb R^n$. In particular, the equation
for $\psi_{J}$ forms an autonomous linear system with
 unique global solution $\psi_{J}(t,u)=\e^{\Bcal^\top_{{J}{J}}
t}\,u_{J}$ for all $u_J\in\mathbb C^n$.
\end{theorem}

Before we prove the theorem, let us illustrate the admissibility
conditions~\eqref{eqadmiss} for the diffusion matrix $\alpha(x)$ for
dimension $d=3$ and the corresponding cases $m =0,1,2,3$. Note that
$\alpha(x)=a+\sum_{i=1}^m x_i\alpha_i$, hence in the case $m=0$ we
have
$$
\alpha(x) \equiv a
$$
for an arbitrary positive semi-definite symmetric $3\times 3$-matrix
$a$. For $m=1$, we have
$$
a = \left(\begin{array}{ccc} 0 & 0 & 0 \\ & + & \ast\\  & & +
\end{array}\right) \, , \quad \alpha_1 = \left(\begin{array}{ccc} +
& \ast & \ast\\ & + & \ast \\  & & + \end{array}\right) \, ,
$$
for $m=2$,
$$
a = \left(\begin{array}{ccc} 0 & 0 & 0 \\ & 0 & 0 \\  & & +
\end{array}\right) \, , \quad \alpha_1=\left(\begin{array}{ccc} + &
0 & \ast \\ & 0 & 0 \\ & & + \end{array}\right) \, , \quad
\alpha_2=\left(\begin{array}{ccc} 0 & 0 & 0 \\ & + & \ast \\ & & +
\end{array}\right) \, ,
$$
and for $m=3$,
$$ a=0\,,\quad
\alpha_1=\left(\begin{array}{ccc} + & 0 &0 \\ & 0 & 0 \\ & & 0
\end{array}\right) \, , \quad \alpha_2 =\left(\begin{array}{ccc} 0 &
0 & 0 \\ & + &0 \\ & & 0 \end{array}\right) \, , \quad
\alpha_3=\left(\begin{array}{ccc} 0 & 0 & 0 \\ & 0 & 0 \\ & & +
\end{array}\right) \, ,
$$
where we leave the lower triangle of symmetric matrices blank, $+$
denotes a non-negative real number and $\ast$ any real number such
that positive semi-definiteness holds.

\begin{proof}[Proof of Theorem~\ref{thmaffine}]
Suppose $X$ is affine. That $a(x)$ and $b(x)$ are of the form
\eqref{eqaffba} follows from Theorem~\ref{thmaffchar}. Obviously,
$a(x)$ is symmetric positive semi-definite for all $x\in\mathbb
R_+^m\times\mathbb R^n$ if and only if $\alpha_j=0$ for all
$j\in{J}$, and $a$ and $\alpha_i$ are symmetric positive
semi-definite for all $i\in{I}$.

We extend the diffusion matrix and drift continuously to $\mathbb
R^d$ by setting
\[ a(x)=a+\sum_{i\in I} x_i^+ \alpha_i\quad\text{and}\quad b(x)=b+\sum_{i\in I} x_i^+
\beta_i+\sum_{j\in J} x_j\beta_j.\]

Now let $x$ be a boundary point of $\mathbb R^m_+\times\mathbb R^n$.
That is, $x_k=0$ for some $k\in{I}$. The stochastic invariance
Lemma~\ref{lemstochinvN} below implies that the diffusion must be
``parallel to the boundary'',
\[ e_k^\top\left(a+\sum_{i\in{I}\setminus\{k\}} x_i\alpha_i \right)e_k=0,\] and the drift
must be ``inward pointing'',
\[ e_k^\top\left(b+\sum_{i\in{I}\setminus\{k\}}x_i\beta_i+\sum_{j\in{J}}x_j\beta_j\right)
\ge 0.\] Since this has to hold for all $x_i\ge 0$,
$i\in{I}\setminus\{k\}$, and $x_j\in\mathbb R$, $j\in{J}$, we obtain
the following set of admissibility conditions
\begin{align*}
a,\,\alpha_i &\;\text{ are symmetric positive semi-definite,}\\
a\,e_k&=0\quad\text{for all $k\in{I}$,}\\
\alpha_i\,e_k&=0\quad\text{for all $i\in{I}\setminus\{k\}$, for all
$k\in{I}$,}\\ \alpha_j&=0\quad\text{for all $j\in{J}$,}\\ b&\in
\mathbb R^m_+\times \mathbb R^n,\\ \beta_i^\top e_k&\ge
0\quad\text{for all $i\in{I}\setminus\{k\}$, for all $k\in{I}$,}\\
\beta_j^\top e_k&=0\quad\text{for all $j\in{J}$, for all $k\in{I}$,}
\end{align*}
which is equivalent to \eqref{eqadmiss}. The form of the
system~\eqref{eqriccfulladmiss} follows by inspection.

Now suppose $a,\,\alpha_i,\,b,\,\beta_i$ satisfy the admissibility
conditions \eqref{eqadmiss}. We show below that there exists a
unique global  solution $(\phi(\cdot,u),\psi(\cdot,u)):\mathbb
R_+\to \mathbb C_-\times\mathbb C^m_-\times\im\mathbb R^n$ of
\eqref{eqriccfulladmiss}, for all $u\in\mathbb C^m_-\times\im\mathbb
R^n$. In particular, $\phi(t,u)+\psi(t,u)^\top x$ has negative real
part for all $t\ge 0$, $u\in\im\mathbb R^d$ and $x\in\mathbb
R^m_+\times\mathbb R^n$. Thus the first part of the theorem follows
from Theorem~\ref{thmaffchar}.

As for the global existence and uniqueness statement, in view of
Lemma~\ref{lemextanalytic1}, it remains to show that $\psi(t,u)$ is
$\mathbb C^m_-\times\im\mathbb R^n$-valued and $t_+(u)=\infty$ for
all $u\in\mathbb C^m_-\times\im\mathbb R^n$. For $i\in I$, denote
the right-hand side of the equation for $\psi_i$ by
\[ R_i(u)= \frac{1}{2}u^\top
    \alpha_i\,u+\beta^\top_{i}u,\]
and observe that
\[ \Re R_i(u) =\frac{1}{2} \Re u^\top \alpha_{i}\,\Re u  -\frac{1}{2}\Im
u^\top\alpha_i\,\Im u + \beta^\top_{i} \Re u . \] Let us denote
$x_I^+=(x_1^+,\dots,x_m^+)^\top$. Since $\Re\psi_J(t,u)=0$, it
follows from the admissibility conditions~\eqref{eqadmiss} and
Corollary~\ref{corlemstochinvN1} below, setting $f(t)=-\Re
\psi(t,u)$,
\[ b_i(t,x)= -\frac{1}{2}\alpha_{i,ii} \left(x_i^+\right)^2
+\frac{1}{2}\Im\psi(t,u)^\top\alpha_i\,\Im\psi(t,u)+\beta_{i,I}^\top
x_I^+,\quad i\in I,\] and $b_j(t,x)=0$ for $j\in J$, that the
solution $\psi(t,u)$ of \eqref{eqriccfulladmiss} has to take values
in $\mathbb C^m_-\times\im\mathbb R^n$ for all initial points $u\in
\mathbb C^m_-\times\im\mathbb R^n$.

Further, for $i\in I$ and $u\in\C^d$, one verifies that
\begin{align*}
 \Re (\overline{u_i} R_i(u)) &=  \frac{1}{2}
\alpha_{i,ii} |u_i|^2\Re u_i +\Re
(\overline{u_i}\,u_i\alpha_{i,{i}{J}}\,u_J)+\frac{1}{2}\Re
(\overline{u_i}\,u_J^\top \alpha_{i,{J}{J}}\,u_J)+\Re
(\overline{u_i}\,\beta^\top_{i} u)\\ &\le \frac{K}{2} \left(1+\|(\Re
u_I)^+\| +\|u_J\|^2\right)\left(1+\|u_I\|^2\right)
\end{align*}
for some finite constant $K$ which does not depend on $u$. We thus
obtain
\begin{align*}
  \partial_t \|\psi_I(t,u)\|^2 &= 2\Re \left(\overline{\psi_I(t,u)}^\top
  R_I\left(\psi_I(t,u),\e^{ \Bcal^\top_{{J}{J}} t}\,u_{J}\right)\right)\\
  &\le  K g(t) \left(1+\|\psi_I(t,u)\|^2\right)
\end{align*}
for
\[ g(t)=\left(1+\|(\Re \psi_I(t,u))^+\|+\|\e^{
\Bcal^\top_{{J}{J}}t}\,u_{J}\|^2\right).\]
 Gronwall's inequality (\cite[(10.5.1.3)]{dieu60}), applied to
$(1+\|\psi_I(t,u)\|^2)$, yields
\begin{equation}\label{extanaeqgron}
\|\psi_I(t,u)\|^2\le\|u_I\|^2  + K\left(1+\|u_I\|^2\right)\int_0^t
g(s) \e^{K\int_s^t g(\xi) \, d\xi }\,ds.
\end{equation}

From above, for all initial points $u\in \mathbb
C^m_-\times\im\mathbb R^n$, we know that $(\Re \psi_I(t,u))^+=0$ and
therefore $t_+(u)=\infty$ by \eqref{extanaeqgron}. Hence the theorem
is proved.
\end{proof}

Now suppose $X$ is affine with characteristics~\eqref{eqaffba}
satisfying the admissibility conditions \eqref{eqadmiss}. In what
follows we show that not only can the functions $\phi(t,u)$ and
$\psi(t,u)$ be extended beyond $u\in\im\mathbb R^d$, but also the
validity of the affine transform formula~\eqref{eqaffdef1} carries
over. This asserts exponential moments of $X(t)$ in particular and
will prove most useful for deriving pricing formulas in affine
factor models.

For any set $U\subset \mathbb R^k$ ($k\in\mathbb N$), we define the
strip
\[ \Scal(U)=\left\{ z\in\mathbb C^k\mid \Re z\in U\right\}\]
in $\mathbb C^k$. The proof of the following theorem builds on
results that are derived in Sections~\ref{secregcharfun} and
\ref{sec: comp} below.

\begin{theorem}\label{thmextanaXXX}
Suppose $X$ is affine with admissible parameters as given in
\eqref{eqadmiss}. Let $\tau>0$. Then
\begin{enumerate}
\item\label{thmextanaXXX1} $\Scal(\Dcal_\mathbb
R(\tau))\subset\Dcal_\mathbb C(\tau)$
\item\label{thmextanaXXX4} $\Dcal_\mathbb R(\tau)=M(\tau)$ where
\[ M(\tau)=\left\{u\in \mathbb R^d\mid\text{$\Ex\left[\e^{u^\top
X^x(\tau)}\right]<\infty$ for all $x\in \R^m_+\times\R^n$}
\right\}.\]
\item\label{thmextanaXXX5} $\Dcal_\mathbb R(\tau)$ and
$\Dcal_\mathbb R$ are convex sets.
\end{enumerate}

Moreover, for all $0\le t\le T$ and $x\in\R^m_+\times\R^n$,
\begin{enumerate}
\setcounter{enumi}{3}
\item\label{thmextanaXXX2} \eqref{eqaffdef1} holds for all $u\in\Scal(\Dcal_\mathbb
R(T-t))$

\item\label{thmextanaXXX3} \eqref{eqaffdef1} holds for all $u\in\mathbb
C^m_-\times\im\mathbb R^n$

\item\label{thmextanaXXX6} $M(t)\supseteq M(T)$.

\end{enumerate}
\end{theorem}

\begin{proof}
We first claim that, for every $u\in\C^d$ with $t_+(u)<\infty$,
there exists some $i\in I$ and some sequence $t_n\uparrow t_+(u)$
such that
\begin{equation}\label{extanaeq2}
  \lim_{n}  (\Re \psi_i(t_n,u))^+=\infty.
\end{equation}
Indeed, otherwise we would have $\sup_{t\in [0,t_+(u))}\|(\Re
\psi_I(t,u))^+\|<\infty$. But then \eqref{extanaeqgron} would imply
$\sup_{t\in [0,t_+(u))} \|\psi_I(t,u)\|<\infty$, which is absurd.
Whence \eqref{extanaeq2} is proved.

In the following, we write
\[ G(u,t,x)=\Ex\left[\e^{u^\top X^x(t)}\right],\quad V(t,x)=\left\{
u\in\mathbb R^d\mid G(u,t,x)<\infty\right\}.\] Since $X$ is affine,
by definition we have $\mathbb R_+\times \im\mathbb
R^d\subset\Dcal_\mathbb C$ and \eqref{eqaffdef1} implies
\begin{equation}\label{eqaffdef1alt}
  G(u,t,x)=\e^{\phi(t,u)+\psi(t,u)^\top x}
\end{equation}
for all $u\in\im\mathbb R^d$, $t\in\mathbb R_+$ and
$x\in\R^m_+\times\R^n$. Moreover, by Lemma~\ref{lemstarshapedEXT}
 , $\Dcal_\R(t)=\Dcal_\C(t)\cap\R^d$ is open and star-shaped
around $0$ in $\R^d$. Hence Lemma~\ref{lemRC3X} implies that
$\Dcal_\R(t)\subset V(t,x)$ and \eqref{eqaffdef1alt} holds for all
$u\in \Dcal_\C(t)\cap\Scal(\Dcal_\R(t))$, for all
$x\in\R^m_+\times\R^n$ and $t\in [0,\tau]$.

Now let $u\in\Dcal_\R(\tau)$ and $v\in \R^d$, and define
\[ \theta^\ast=\inf\{\theta\in\R_+\mid u+\im \theta
v\notin\Dcal_\C(\tau)\}.\] We claim that $\theta^\ast=\infty$.
Arguing by contradiction, assume that $\theta^\ast<\infty$. Since
$\Dcal_\C(\tau)$ is open, this implies $u+\im \theta^\ast v\notin
\Dcal_\C(\tau)$, and thus \begin{equation}\label{thmextanaXXX1cont}
t_+(u+\im \theta^\ast v)\le\tau.
\end{equation}
On the other hand, since $\Dcal_\R(\tau)$ is open,
$(1+\epsilon)u\in\Dcal_\R(\tau)$ for some $\epsilon>0$. Hence
\eqref{eqaffdef1alt} holds and $G(t,(1+\epsilon)u,x)$ is uniformly
bounded in $t\in [0,\tau]$, by continuity of $\phi(t,(1+\epsilon)u)$
and $\psi(t,(1+\epsilon)u)$ in $t$. We infer that the class of
random variables $\{\e^{(u+\im \theta^\ast v)^\top X(t)}\mid t\in
[0,\tau]\}$ is uniformly integrable, see \cite[13.3]{will91}. Since
$X(t)$ is continuous in $t$,  we conclude by Lebesgue's convergence
theorem that $G(t,u+\im \theta^\ast v,x)$ is continuous in $t\in
[0,\tau]$, for all $x\in\R^m_+\times \R^n$. But for all $t<t_+(u+\im
\theta^\ast v)$ we have $(t,u+\im \theta^\ast v)\in
\Dcal_\C(t)\cap\Scal(\Dcal_\R(t))$, and thus \eqref{eqaffdef1alt}
holds for all $x\in\R^m_+\times \R^n$. In view of \eqref{extanaeq2},
this contradicts \eqref{thmextanaXXX1cont}. Whence
$\theta^\ast=\infty$ and thus $u+\im v\in\Dcal_\C(\tau)$. This
proves \ref{thmextanaXXX1}\footnote{For an alternative proof of the
above, see remark \ref{alternative proof}}.

Applying the above arguments to\footnote{Here we use the Markov
property of $X$, see \cite[Theorem 5.4.20]{kar/shr:91}.}
$\Ex\left[\e^{u^\top X(T)}\mid\Fcal_t\right]=G(T-t,u,X(t))$ with
$T=t+\tau$ yields \ref{thmextanaXXX2}. Part~\ref{thmextanaXXX3}
follows, since, by Theorem~\ref{thmaffine},
$\C^m_-\times\im\R^n\subset\Scal(\Dcal_\R(t))$ for all $t\in\R_+$.

As for \ref{thmextanaXXX4}, we first let $u\in \mathcal D_{\mathbb
R}(\tau)$. From part~\ref{thmextanaXXX2} it follows that $u\in
M(\tau)$. Conversely, let $u\in M(\tau)$, and define
$\theta^\ast=\sup\{\theta\ge 0\mid \theta u\in \Dcal_\R(\tau)\}$. We
have to show that $\theta^\ast>1$. Assume, by contradiction, that
$\theta^\ast\le 1$. From Lemma~\ref{lemstarshapedEXT}  , we know
that there exists some $x^\ast\in\R^m_+\times\R^n$ such that
\begin{equation}\label{eqthmextanaXXX1}
  \lim_{\theta\uparrow\theta^\ast} \phi(\tau,\theta u)+\psi(\tau,\theta
  u)^\top x^\ast=\infty.
\end{equation}
On the other hand, from part~\ref{thmextanaXXX2} and Jensen's
inequality, we obtain
\[ \e^{\phi(\tau,\theta  u)+\psi(\tau,\theta  u)^\top
x^\ast}=G(\tau,\theta u,x^\ast)\le G(\tau,u,x^\ast)^\theta\le
G(\tau,u,x^\ast)<\infty\] for all $\theta<\theta^\ast$. But this
contradicts \eqref{eqthmextanaXXX1}, hence $u\in\Dcal_\R(\tau)$, and
part~\ref{thmextanaXXX4} is proved. Since $M(\tau)$ is convex, this
also implies \ref{thmextanaXXX5}. Finally, part~\ref{thmextanaXXX6}
follows from part~\ref{thmextanaXXX4} and
Lemma~\ref{lemextanalytic1}. Whence the theorem is proved.
\end{proof}

\begin{remark}\label{remglasskim}
Glasserman and Kim~\cite{Glasserman} proved the equality in
Theorem~\ref{thmextanaXXX}~\ref{thmextanaXXX4}, and the validity of
the transform formula \eqref{eqaffdef1} for all $u$ in an open
neighborhood of $\Dcal_\mathbb R(T-t)$ in $\C^d$, under the
additional assumption that $\Bcal$ has strictly negative
eigenvalues. That assumption, however, excludes the simple Heston
stochastic volatility model in Section~\ref{secheston} below.
\end{remark}

\begin{remark}
In Keller-Ressel~\cite[Theorem 3.18 and Lemma 3.19]{keller} it is
shown that
\[ M(\tau+\epsilon)\subseteq \Dcal_\R(\tau) \]
for all $\epsilon>0$, for a more general class of affine Markov
processes $X^x$. Obviously, in our framework, this is implied by
parts~\ref{thmextanaXXX4} and \ref{thmextanaXXX6} of
Theorem~\ref{thmextanaXXX}.
\end{remark}

\begin{remark}
The convexity property of the maximal domain stated in
Theorem~\ref{thmextanaXXX}~\ref{thmextanaXXX5} represents a
non-trivial result for ordinary differential equations. Only in the
mid 1990s have corresponding convexity results been derived in the
analysis literature, see Lakshmikantham et al.\ \cite{lakshawal96}.
\end{remark}

\section{Discounting and Pricing in Affine Models}\label{secdiscpri}

We let $X$ be affine on the canonical state space $\mathbb
R^m_+\times \mathbb R^n$ with admissible parameters
$a,\alpha_i,b,\beta_i $ as given in \eqref{eqadmiss}. Since we are
interested in pricing, and to avoid a change of measure, we
interpret $\Pa=\Q$ as risk-neutral measure in what follows.

A short rate model of the form
\begin{equation} \label{r}
r(t) = c + \gamma^\top X(t),
\end{equation}
for some constant parameters $c \in \mathbb{R}$ and $\gamma \in
\mathbb{R}^d$, is called an affine short rate model. Special cases,
for dimension $d=1$, are the Vasi\v{c}ek and Cox--Ingersoll--Ross
short rate models. We recall that an affine term structure model
always induces an affine short rate model.

Now consider a $T$-claim with payoff $f(X(T))$. Here $f:\mathbb
R^m_+\times\mathbb R^n\to \mathbb R$ denotes a measurable payoff
function, such that $f(X(T))$ meets the required integrability
conditions
\[ \Ex\left[ \e^{-\int_0^T r(s)\,ds}
\,|f(X(T))|\right]<\infty.\] Its arbitrage price at time $t\le T$ is
then given by
\begin{equation}\label{fPriceqn}
\pi(t)=\Ex\left[ \e^{-\int_t^T r(s)\,ds}
\,f(X(T))\mid\Fcal_t\right].
\end{equation}
A particular example is the $T$-bond with $f\equiv 1$. Our aim is to
derive an analytic, or at least numerically tractable, pricing
formula for \eqref{fPriceqn}. To this end we shall make use of a
change of numeraire technique to price, e.g., Bond options and
caplets. Denote the risk free bank account by $B(t):=\e^{\int_0^t
r(s)ds}$. For fixed $T>0$ it is easily observed that
\[
\frac{1}{P(0,T)B(T)}>0\quad\textit{and}\quad\mathbb
E\left[\frac{1}{P(0,T)B(T)}\right]=1,
\]
hence we may introduce an equivalent probability measure $\mathbb
Q^T\sim \mathbb Q$ on $\mathcal F_T$ by its Radon-Nikodym derivative
\[
\frac{d\mathbb Q^T}{d\mathbb Q}=\frac{1}{P(0,T)B(T)}.
\]
$\mathbb Q^T$ is called the $T$-forward measure. Note that for
$t\leq T$,
\begin{equation}\label{eq: conditioned radon nikodym der}
\left.\frac{d\mathbb Q^T}{d\mathbb Q}\right|_{\mathcal F_t}=\mathbb
E[\left.\left[\frac{1}{P(0,T)B(T)}\right| \mathcal
F_t\right]=\frac{P(t,T)}{P(0,T)B(t)}.
\end{equation}
As a first step towards establishing useful pricing formulas, we
derive a formula for the $\Fcal_t$-conditional characteristic
function of $X(T)$ under $\mathbb Q^T$, which up to normalization
with $\Ex\left[ \e^{-\int_t^T r(s)\,ds} \mid\Fcal_t\right]$ equals,
\begin{equation}\label{eqpriceaffinepre}
 \Ex\left[ \e^{-\int_t^T r(s)\,ds} \,\e^{u^\top
X(T)}\mid\Fcal_t\right],\quad u\in\im\mathbb R^d
\end{equation}
(use equation\ \eqref{eq: conditioned radon nikodym der}.

Note that the following integrability
condition~{\ref{thmaffdisc1aX}} is satisfied in particular if $r$ is
uniformly bounded from below, that is, if $\gamma\in\mathbb
R^m_+\times\{0\}$.
\begin{theorem}\label{thmaffdisc1}
Let $\tau>0$. The following statements are equivalent:
\begin{enumerate}
\item\label{thmaffdisc1aX} $\Ex\left[ \e^{-\int_0^\tau
r(s)\,ds} \right]<\infty$ for all $x\in\mathbb R^m_+\times\mathbb
R^n$.

\item\label{thmaffdisc1bX}  There exists a unique solution
$(\Phi(\cdot,u),\Psi(\cdot,u)):[0,\tau]\to \mathbb C\times\mathbb
C^d$ of
\begin{equation}\label{canricceqexXX}
\begin{aligned}
\partial_t\Phi(t,u)&= \frac{1}{2} \Psi_J(t,u)^\top
a_{JJ}\,\Psi_J(t,u)+b^\top\Psi(t,u)-c,\\ \Phi(0,u)&=0,\\
\partial_t\Psi_i(t,u)&=\frac{1}{2} \Psi (t,u)^\top
\alpha_{i }\,\Psi (t,u) +\beta_i^\top\Psi(t,u)-\gamma_i  ,\quad i\in
I,\\
\partial_t\psi_{J}(t,u)&=\Bcal_{{J}{J}}^\top
  \Psi_{J}(t,u)-\gamma_J ,\\
  \Psi(0,u)&=u
\end{aligned}
\end{equation}
for $u=0$.

\end{enumerate}

In either case, there exists an open convex neighborhood $U$ of $0$
in $\R^d$ such that the system of Riccati
equations~\ref{canricceqexXX} admits a unique solution
$(\Phi(\cdot,u),\Psi(\cdot,u)):[0,\tau]\to \mathbb C\times\mathbb
C^d$ for all $u\in \Scal(U)$, and \eqref{eqpriceaffinepre} allows
the following affine representation
\begin{equation}\label{eqpricecharXX}
 \Ex\left[ \e^{-\int_t^T r(s)\,ds}
\,\e^{u^\top
X(T)}\mid\Fcal_t\right]=\e^{\Phi(T-t,u)+\Psi(T-t,u)^\top X(t)}
\end{equation}
for all $u\in\Scal(U)$, $t\le  T\le t+ \tau$ and $x\in\mathbb
R^m_+\times\mathbb R^n$.
\end{theorem}

\begin{proof}
We first enlarge the state space and consider the real-valued
process
\[ Y(t)=y+\int_0^t\left(c+\gamma^\top X(s)\right)ds,\quad
y\in\mathbb R.\] A moment's reflection reveals that
$X'=\left(\begin{array}{c}X\\ Y\end{array}\right)$ is an $\mathbb
R^m_+\times \mathbb R^{n+1}$-valued diffusion process with diffusion
matrix $a'+\sum_{i\in I} x_i\alpha_i'$ and drift $b'+\Bcal' x'$
where
\[a'=\left(\begin{array}{cc} a & 0 \\
0 & 0\end{array}\right),\quad \alpha_i'=\left(\begin{array}{cc}
\alpha_i & 0 \\ 0 & 0\end{array}\right),\quad
b'=\left(\begin{array}{c} b \\ c\end{array}\right),\quad
\Bcal'=\left(\begin{array}{cc} \Bcal & 0 \\ \gamma^T &
0\end{array}\right) \] form admissible parameters. We claim that
$X'$ is an affine process.

Indeed, the candidate system of Riccati equations reads
\begin{equation}\label{canricceqexXXX}
\begin{aligned}
\partial_t\phi'(t,u,v)&= \frac{1}{2} \psi'_J(t,u,v)^\top
a_{JJ}\,\psi'_J(t,u,v)+b^\top\psi'_{\{1,\dots,d\}}(t,u,v)+\text{\fbox{$cv$}},\\
\phi'(0,u,v)&=0,\\
\partial_t\psi'_i(t,u,v)&=\frac{1}{2} \psi' (t,u,v)^\top
\alpha_{i }\,\psi' (t,u,v)
+\beta_i^\top\psi'(t,u,v)+\text{\fbox{$\gamma_i v$}},\quad i\in I,\\
\partial_t\psi'_{J}(t,u,v)&=\Bcal_{{J}{J}}^\top
  \psi'_{J}(t,u,v)+\text{\fbox{$\gamma_J v$}},\\
  \partial_t\psi'_{d+1}(t,u,v)&=0,\\
  \psi'(0,u,v)&=\left(\begin{array}{c} u\\ v\end{array}\right).
\end{aligned}
\end{equation}
Here we replaced the constant solution $\psi'_{d+1}(\cdot,u,v)\equiv
v$ by $v$ in the boxes. Theorem~\ref{thmaffine} carries over and
asserts a unique global $\mathbb C_-\times\mathbb
C^m_-\times\im\mathbb R^{n+1}$-valued solution
$(\phi'(\cdot,u,v),\psi'(\cdot,u,v))$ of \eqref{canricceqexXX} for
all $ (u,v) \in\mathbb C^m_-\times\im\mathbb R^{n}\times\im\mathbb
R$. The second part of Theorem~\ref{thmaffchar} thus asserts that
$X'$ is affine with conditional characteristic function
\[\Ex\left[ \e^{u^\top  X(T)+v Y(T)}\mid\Fcal_t\right] =
\e^{\phi'(T-t,u,v)+\psi'(T-t,u,v)^\top   X (t)+vY(t)}\] for all $
(u,v) \in\mathbb C^m_-\times\im\mathbb R^{n}\times\im\mathbb R$ and
$t\le T$.

The theorem now follows from Theorem~\ref{thmextanaXXX} once we set
$\Phi(t,u)=\phi'(t,u,-1)$ and
$\Psi(t,u)=\psi'_{\{1,\dots,d\}}(t,u,-1)$.
\end{proof}

Suppose, for the rest of this section, that either condition
{\ref{thmaffdisc1aX}} or {\ref{thmaffdisc1bX}} of
Theorem~\ref{thmaffdisc1} is met. As immediate consequence of
Theorem~\ref{thmaffdisc1}, we obtain the following explicit price
formulas for $T$-bonds in terms of $\Phi$ and $\Psi$.

\begin{corollary}\label{charfunctSforwardmeasure}
For any maturity $T\le \tau$, the $T$-bond price at $t\le T$ is
given as
\[ P(t,T)=\e^{-A(T-t)-B(T-t)^\top X(t)} \]
where we denote
\[ A(t)=-\Phi(t,0),\quad B(t)=-\Psi(t,0).\]

Moreover, for $t\le T\le S\le\tau$, the $\Fcal_t$-conditional
characteristic function of $X(T)$ under the $S$-forward measure
$\Q^S$ is given by
\begin{equation}\label{chfctSTforBB}
\Ex_{\Q^S}\left[\e^{u^\top X(T)}\mid\Fcal_t\right]
=\frac{\e^{-A(S-T)+\Phi(T-t,u-B(S-T))+\Psi(T-t,u-B(S-T))^\top
X(t)}}{P(t,S)}
\end{equation}
for all $u\in \Scal(U+B(S-T))$, where $U$ is the neighborhood of $0$
in $\mathbb R^d$ from Theorem~\ref{thmaffdisc1}.
\end{corollary}

\begin{proof}
The bond price formula follows from \eqref{eqpricecharXX} with
$u=0$.

Now let $t\le T\le S\le\tau$ and $u\in \Scal(U+B(S-T))$. We obtain
from \eqref{eqpricecharXX} by nested conditional expectation
\begin{align*}
\Ex\left[\e^{-\int_t^S r(s)\,ds}\e^{u^\top
X(T)}\mid\Fcal_t\right]&=\Ex\left[\e^{-\int_t^T
r(s)\,ds}\,\Ex\left[\e^{-\int_T^S
r(s)\,ds}\mid\Fcal_T\right]\e^{u^\top X(T)}\mid\Fcal_t\right]\\
&=\frac{\Ex\left[\e^{-\int_t^T r(s)\,ds}\e^{ (u-B(S-T))^\top
X(T)}\mid\Fcal_t\right]}{\e^{A(S-T)}}\\
&=\frac{\e^{\Phi(T-t,u-B(S-T))+\Psi(T-t,u-B(S-T))^\top
X(t)}}{\e^{A(S-T)}}.
\end{align*}
Normalizing by $P(t,S)$ yields \eqref{chfctSTforBB}.
\end{proof}

For more general payoff functions $f$, we can proceed as follows.

\begin{itemize}

\item Either we recognize the $\Fcal_t$-conditional distribution, say $q(t,T,dx)$, of
$X(T)$ under the $T$-forward measure from its characteristic
function~\eqref{chfctSTforBB}. Or we derive $q(t,T,dx)$ via
numerical inversion of the characteristic
function~\eqref{chfctSTforBB}, using e.g.\ fast Fourier transform
(FFT). Then compute the price \eqref{fPriceqn} by integration of $f$
\begin{equation}\label{eqpiaffcompX}
\pi(t) = P(t,T)\,\int_{\mathbb R^m_+\times\mathbb R^n}
f(x)\,q(t,T,dx).
\end{equation}
Examples are given in Section~\ref{secbondopaffine} below.

\item Or suppose $f$ can be expressed by
\begin{equation} \label{eqfinvXX}
 f(x) = \int_{\mathbb R^d} \e^{(u+\im y)^\top x} \,\widetilde{f}(y)\,dy
\end{equation}
for some integrable function $\widetilde{f}:\mathbb R^d\to\mathbb C$
and some constant $u\in U$. Then we may apply Fubini's theorem to
change the order of integration, which gives
\begin{equation}\label{exopthestonXX}
\begin{aligned}
  \pi(t)&=\Ex\left[ \e^{-\int_t^T r(s)\,ds}
\,\int_{\mathbb R^d} \e^{(u+\im y)^\top X(T)}
\widetilde{f}(y)\,dy\mid\Fcal_t\right]\\ &=\int_{\mathbb R^d}
\Ex\left[ \e^{-\int_t^T r(s)\,ds} \, \e^{(u+\im y)^\top X(T)}
\mid\Fcal_t\right] \widetilde{f}(y)\,dy\\ &= \int_{\mathbb R^d}
\e^{\Phi(T-t,u+\im y)+\Psi(T-t,u+\im y)^\top X(t)}
\,\widetilde{f}(y)\,dy.
\end{aligned}
\end{equation}
This integral can be numerically computed. An example is given in
Section~\ref{secheston} below.
\end{itemize}

The function $\widetilde{f}$ in \eqref{eqfinvXX} can be found by
Fourier transformation, as the following classical result indicates.
\begin{lemma}\label{lemfourier}
Let $f:\mathbb R^d\to\mathbb C$ be a measurable function and
$u\in\mathbb R^d$ be such that the function $h(x)=e^{-u^\top
x}\,f(x)$ and its Fourier transform
\[ \hat{h}(y) =  \int_{\mathbb R^d} h(x)\,\e^{-\im y^\top
x}\,dx \] are integrable on $\mathbb R^d$. Then \eqref{eqfinvXX}
holds for almost all $x\in\mathbb R^d$ for
\[ \widetilde{f}=\frac{1}{(2\pi)^d} \hat{h} .\]

Moreover, the right hand side of \eqref{eqfinvXX} is continuous in
$x$. Hence, if $f$ is continuous then \eqref{eqfinvXX} holds for all
$x\in\mathbb R^d$.
\end{lemma}

\begin{proof}
From Fourier analysis, see \cite[Chapter I, Corollary
1.21]{stewei71}, we know that
\[ h(x)= \frac{1}{(2\pi)^d} \int_{\mathbb R^d} \e^{ \im y^\top x}
\,\hat{h}(y)\,dy \] for almost all $x\in\mathbb R^d$. Multiplying
both sides with $\e^{u^\top x}$ yields the first claim.

From the Riemann--Lebesgue Theorem (\cite[Chapter I, Theorem
1.2]{stewei71}) we know that the right hand side of \eqref{eqfinvXX}
is continuous in $x$.
\end{proof}

An example is the continuous payoff function
\[f(x)=(\e^x-K)^+\] of a European call option with strike price $K$
on the underlying stock price $\e^{L }$, where $L$ may be any affine
function of $X$. Fix a real constant ${p}>1$. Then
$h(x)=\e^{-{p}x}f(x)$ is integrable on $\mathbb R$. An easy
calculation shows that its Fourier transform
\[ \hat{h}(y) =\int_\mathbb R \e^{-{p}x}f(x)\, \e^{-\im y x}\,dx = \frac{K^{1-{p}-\im
y}}{({p}+\im y)({p}+\im y-1)} \] is also integrable on $\mathbb R$.
In view of Lemma~\ref{lemfourier}, we thus conclude that, for
${p}>1$,
\begin{equation}\label{exoptheston}
 (\e^x-K)^+ = \frac{1}{ 2\pi } \int_\mathbb R \e^{({p}+\im y) x}\,\frac{K^{1-{p}-\im
y}}{({p}+\im y)({p}+\im y-1)}\, dy,
\end{equation}
which is of the desired form~\eqref{eqfinvXX}. We will apply this
for the Heston stochastic volatility model in
Section~\ref{secheston} below.

A related example is the following
\begin{equation}\label{eqrelexamf}
(\e^x-K)^+ -\e^x =\frac{1}{ 2\pi } \int_\mathbb R \e^{({p}+\im y)
x}\,\frac{K^{1-{p}-\im y}}{({p}+\im y)({p}+\im y-1)}\, dy,
\end{equation}
which holds for all $0<{p}<1$.

More examples of payoff functions with integral representation,
including the above, can be found in \cite{hubkalkra06}.

\section{Bond Option Pricing in Affine Models}\label{secbondopaffine}

We can further simplify formula \eqref{eqpiaffcompX} for a European
call option on a $S$-bond with expiry date $T<S$ and strike price
$K$. The payoff function is
\[ f(x) = \left(\e^{-A(S-T)-B(S-T)^\top
x}-K\right)^+.\] We can decompose \eqref{fPriceqn},
\begin{equation}\label{eqpidecSTfmC}
 \pi^C(t;T,S)= P(t,S) \Q^S[E \mid\Fcal_t]
 -K P(t,T)\Q^T[E \mid\Fcal_t]
\end{equation}
for the event $E=\{B(S-T)^\top X(T)\le -A(S-T)-\log K \}$. The
pricing of this bond option boils down to the computation of the
probability of the event $E$ under the $S$- and $T$-forward
measures.

Similarly, the value of a put equals
\begin{equation}\label{eqpidecSTfmP}
 \pi^P(t; T,S)= K P(t,T)\Q^T[E^c \mid\Fcal_t]-P(t,S) \Q^S[E^c \mid\Fcal_t]
\end{equation}
for the event $E^c=\Omega\setminus E=\{B(S-T)^\top X(T)>
-A(S-T)-\log K \}$.

In the following two subsections, we illustrate this approach for
the Vasi\v{c}ek and Cox--Ingersoll--Ross short rate models.

\subsection{Example: Vasi\v{c}ek Short Rate
Model}\label{subsecvasiaff}

The state space is $\mathbb R$, and we set $r=X$ for the Vasi\v{c}ek
short rate model \[dr=(b+\beta r)\,dt +\sigma\,dW.\] The
system~\eqref{canricceqexXX} reads
\begin{align*}
 \Phi(t,u)&= \frac{1}{2} \sigma^2 \int_0^t \Psi^2(s,u)\,ds +b \int_0^t
\Psi(s,u)\,ds\\
\partial_t\Psi (t,u)&=\beta \Psi (t,u)-1,\notag\\
  \Psi(0,u)&=u\notag
\end{align*}
which admits a unique global solution with
\begin{align*}
\Psi(t,u)&=\e^{\beta t}u - \frac{\e^{\beta t}-1}{\beta}\\
\Phi(t,u)&= \frac{1}{2}\sigma^2\left(\frac{u^2}{2\beta }(\e^{2\beta
t}-1)+\frac{1}{2\beta^3}(\e^{2\beta t}-4\e^{\beta t}+2\beta
t+3)\right.\\-&\left.\frac{u}{\beta^2}(\e^{2\beta t}-2\e^{\beta
t}+2\beta)\right)+b\left(\frac{\e^{\beta
t}-1}{\beta}\,u+\frac{\e^{\beta t}-1-\beta t}{\beta^2}\right)
\end{align*}
for all $u\in\mathbb C$. Hence \eqref{eqpricecharXX} holds for all
$u\in\mathbb C$ and $t\le T$. In particular, by Corollary
\ref{charfunctSforwardmeasure}, the bond prices $P(t,T)$ can be
determined by $A$ and $B$,
\begin{align*}
B(t)&=-\Psi(t,0)=\frac{\e^{\beta t}-1}{\beta},\\
A(t)&=-\Phi(t,0)=-\frac{\sigma^2}{4\beta^3}(\e^{2\beta t}-4\e^{\beta
t}+2\beta t+3)+b\frac{\e^{\beta t}-1-\beta t}{\beta^2}.
\end{align*}
Hence, under the $S$-forward measure, $r(T)$ is
$\Fcal_t$-conditionally Gaussian distributed with (cf.~
\cite{brimer:06}, chapter 3.2.1)
\begin{align*}
\mathbb E_{\mathbb Q^S}[r(T)\mid\mathcal
F_t]&=r(t)e^{-\beta(T-s)}+M^S(t,T),
 \\\Var(r(T)\mid\mathcal F_t)&=\sigma^2\frac{\e^{2\beta
(T-t)}-1}{2\beta},
\end{align*}
where $M^S$ is defined by
\[
M^S(t,T)=(\frac{b}{\beta}
-\frac{\sigma^2}{2\beta^2})(1-\e^{-\beta(T-t)})
+\frac{\sigma^2}{2\beta^2}\left[\e^{-\beta(S-T)}-\e^{-\beta(S+T-2t)}\right].
\]
The bond option price formula for the Vasi\v{c}ek short rate model
can now be derived via \eqref{eqpidecSTfmC} and
\eqref{eqpidecSTfmP}.

\subsection{Example: Cox--Ingersoll--Ross Short Rate
Model}\label{subsecciraff}

 The state space
is $\mathbb R_+$, and we set $r=X$ for the Cox--Ingersoll--Ross
short rate model \[dr=(b+\beta r)\,dt +\sigma\sqrt{r}\,dW.\] The
system~\eqref{canricceqexXX} reads
\begin{equation}\label{ricceqcirUU}
\begin{aligned}
 \Phi(t,u)&= b \int_0^t \Psi(s,u)\,ds,\\
\partial_t\Psi (t,u)&=\frac{1}{2} \sigma^2   \Psi^2(t,u)+\beta \Psi (t,u)-1,\\
  \Psi(0,u)&=u.
\end{aligned}
\end{equation}
By Lemma~\ref{lemsolricc} below, there exists a unique solution
$(\Phi(\cdot,u),\Psi(\cdot,u)):\mathbb R_+\to \mathbb
C_-\times\mathbb C_-$, and thus \eqref{eqpricecharXX} holds, for all
$u\in\mathbb C_-$ and $t\le T$. The solution is given explicitly as
\begin{align*}
  \Phi(t,u)&=\frac{2b}{\sigma^2}\log\left(\frac{ L_5(t)
  }{L_3(t)-L_4(t) u  }\right)\\
  \Psi(t,u)&=-\frac{L_1(t) - L_2(t) u }
  {L_3(t)-L_4(t) u }
\end{align*}
where ${\lambda}=\sqrt{\beta^2+2\sigma^2}$ and
\begin{align*}
  L_1(t)&=2 \left(\e^{{\lambda}
  t}-1\right)\\
  L_2(t)&= {\lambda}\left(\e^{{\lambda} t}+1\right)+\beta\left(\e^{{\lambda}
  t}-1\right) \\
  L_3(t)&={\lambda}\left(\e^{{\lambda} t}+1\right)-\beta\left(\e^{{\lambda} t}-1\right) \\
  L_4(t)&=\sigma^2\left(\e^{{\lambda}
  t}-1\right)\\
  L_5(t)&=2{\lambda}\e^{\frac{({\lambda}-\beta)t}{2}}.
\end{align*}

Some tedious but elementary algebraic manipulations show that the
$\Fcal_t$-conditional characteristic function of $r(T)$ under the
$S$-forward measure $\Q^S$ is given by
\[\Ex_{\Q^S}\left[\e^{u r(T)}\mid\Fcal_t\right]
=\frac{\e^{-C_2(t,T,S)r(t)+\frac{C_2(t,T,S)r(t)}{1-C_1(t,T,S)u}}}{(1-C_1(t,T,S)u)^\frac{2b}{\sigma^2}}\]
where
\[C_1(t,T,S) =\frac{L_3(S-T)L_4(T-t)}{2{\lambda} L_3(S-t)},\quad
  C_2(t,T,S)= \frac{L_2(T-t)}{L_4(T-t)}-\frac{L_1(S-t)}{L_3(S-t)}\]

Comparing this with Lemma~\ref{lemnoncentrchi} below, we conclude
that the $\Fcal_t$-conditional distribution of the random variable
$2 r(T)/C_1(t,T,S)$ under the $S$-forward measure $\Q^S$ is
noncentral $\chi^2$ with $\frac{4b}{\sigma^2}$ degrees of freedom
and parameter of noncentrality $2 C_2(t,T,S)r(t)$. Combining this
with \eqref{eqpidecSTfmC}--\eqref{eqpidecSTfmP}, we obtain explicit
European bond option price formulas.

As an application, we now compute cap prices. Let us consider a cap
with strike rate $\kappa$ and tenor structure
$1/4=T_0<T_1<\dots<T_n$, with $T_i-{T_{i-1}}=1/4$. Here, as usual,
$T_i$ denote the settlement dates and $T_{i-1}$ the reset dates for
the $i$th caplet, $i=1,\dots,n$ and $T_n$ is the maturity of the
cap. It is well known that the cash flow of a $i$th caplet at time
$T_i$ equals the $(1+\kappa/4)$ multiple of the cash-flow at
$T_{i-1}$ of a put option on the $T_i$-bond with strike price
$1/(1+\kappa/4)$. Hence the cap price equals
\[
Cp=\sum_{i=1}^n Cpl(i)=(1+  \kappa/4)\sum_{i=1}^n
P(0,T_{i-1})\mathbb E_{\mathbb
Q^{T_{i-1}}}\left[\left(\frac{1}{1+\kappa/4}-P(T_{i-1},T_i)\right)^+\right].
\]

In practice, cap prices are often quoted in Black implied
volatilities. By definition, the implied volatility $\sigma_B>0$ is
the number, which, plugged into Black's formula, yields the cap
value $Cp=\sum_{i=1}^n Cpl(i)$, where the $i$th caplet price is
given as
\[
Cpl(i)=\frac{1}{4}
P(0,T_i)(F(T_{i-1},T_i)\Phi(d_1(i))-\kappa\Phi(d_2(i)))
\]
with
\[
d_{1,2}(i)=\frac{\log\left(\frac{F(T_{i-1},T_i)}{\kappa}\right)\pm\frac{\sigma_B^2}{2}(T_{i-1}-t)}{\sigma_B\sqrt{T_{i-1}-t}}.
\]
where
$F(T_{i-1},T_i)=4\left(\frac{P(0,T_{i-1})}{P(0,T_{i})}-1\right)$
denotes the corresponding simple forward rate.

As parameters for the CIR model we assume
\[
\sigma^2=0.033,\quad b=0.08,\quad \beta=-0.9,\quad r_0=0.08.
\]
In Table~\ref{tabcirprices} we summarize the ATM\footnote{The cap
with maturity $T_n$ is at-the-money (ATM) if its strike rate
$\kappa$ equals the prevailing forward swap rate
$4(P(0,T_0)-P(0,T_n))/\sum_{i=1}^n P(0,T_i)$.} cap prices and
implied volatilities for various maturities.

\begin{table}[h!]
\caption{ATM cap prices for the CIR model}\label{tabcirprices}
\begin{center}
\begin{tabular}{c*{3}{|l}}
                \textit{Maturity Years}& \textit{strike rate}&
                \textit{cap
                price}&\textit{implied volatility}\\\hline
             1    &          0.0843    &    0.0073    &       0.4506       \\\hline
              2    &          0.0855    &    0.0190    &       0.3720       \\\hline
              3    &          0.0862    &    0.0302    &       0.3226       \\\hline
              4    &          0.0866    &    0.0406    &       0.2890        \\\hline
              5    &          0.0868    &    0.0501    &       0.2647        \\\hline
              6    &          0.0870    &    0.0588    &       0.2462        \\\hline
              7    &          0.0871    &    0.0668    &       0.2316        \\\hline
              8    &          0.0872    &    0.0742    &       0.2198       \\\hline
              9    &          0.0873    &    0.0809    &       0.2100        \\\hline
             10    &          0.0873    &    0.0871    &       0.2017       \\\hline
             15    &          0.0875   &    0.1110    &       0.1744       \\\hline
             20    &          0.0876    &    0.1265    &       0.1594       \\\hline
             25    &          0.0876    &    0.1365    &       0.1502        \\\hline
             30    &          0.0876    &    0.1430    &       0.1442       \\\hline

\end{tabular}
\end{center}
\end{table}

\begin{lemma}[Noncentral $\chi^2$-Distribution]\label{lemnoncentrchi}
The noncentral $\chi^2$-distribution with $\delta>0$ degrees of
freedom and noncentrality parameter $\zeta>0$ has density function
\[ f_{\chi^2(\delta,\zeta)}(x)= \frac{1}{2} \e^{-\frac{x+\zeta}{2}}
\left(\frac{x}{\zeta}\right)^{\frac{\delta}{4}-\frac{1}{2}}\,I_{\frac{\delta}{2}-1}(\sqrt{\zeta
x}),\quad x\ge 0\]  and characteristic function
\[ \int_{\mathbb R_+} \e^{ux}\,f_{\chi^2(\delta,\zeta)}(x)\,dx
= \frac{\e^{\frac{\zeta u}{1-2u}}}{(1-2u)^\frac{\delta}{2}},\quad
u\in\mathbb C_-.\] Here $I_\nu(x)=\sum_{j\ge
0}\frac{1}{j!\Gamma(j+\nu+1)}\left(\frac{x}{2}\right)^{2j+\nu}$
denotes the modified Bessel function of the first kind of order
$\nu>-1$.
\end{lemma}
\begin{proof}
See e.g.\ \cite{Johnson}.
\end{proof}

\begin{lemma}\label{lemsolricc}
Consider the Riccati differential equation
\begin{equation}\label{eq: scalar Riccati}
\partial_t G=A G^2 + B G - C,\quad G(0,u)=u,
\end{equation}
where $A,B,C\in\mathbb C$ and $u\in\C$, with $A\neq 0$ and
$B^2+4AC\in\mathbb C\setminus \mathbb R_-$. Let $\sqrt{\cdot}$
denote the analytic extension of the real square root to $\mathbb
C\setminus \mathbb R_-$, and define ${\lambda}=\sqrt{B^2+4AC}$.

\begin{enumerate}
\item \label{item1: lemsolricc} The function
\begin{equation}\label{sol Ric: form1}
G(t,u)=-\frac{2 C \left(e^{{\lambda} t}-1\right) -
\left({\lambda}\left(e^{{\lambda} t}+1\right)+B\left(e^{{\lambda}
t}-1\right)\right)u}{{\lambda}\left(e^{{\lambda}
t}+1\right)-B\left(e^{{\lambda} t}-1\right)-2A\left(e^{{\lambda}
t}-1\right)u}
\end{equation}
is the unique solution of equation \eqref{eq: scalar Riccati} on its
maximal interval of existence $[0,t_+(u))$. Moreover,
\begin{equation}\label{eq: integral of G}
\int_0^t G(s,u){d}s=\frac{1}{A}\log\left(\frac{2\lambda
e^{\frac{\lambda-B}{2}t}}{\lambda (e^{\lambda t}+1)-B(e^{\lambda
t}-1)-2A(e^{\lambda t}-1) u}\right).
\end{equation}
\item \label{item2: lemsolricc}  If, moreover, $A> 0$, $B\in\mathbb R$, $\Re(C)\ge 0$ and
$u\in\mathbb C_-$ then $t_+(u)=\infty$ and $G(t,u)$ is
$\C_-$-valued.
\end{enumerate}
\end{lemma}
\begin{proof}
\ref{item1: lemsolricc}: Recall that the square root
$\sqrt{z}:=e^{1/2\log(z)}$ is the well defined analytic extension of
the real square root to $\mathbb C\setminus \mathbb R_-$, through
the main branch of the logarithm which can be written in the form
$\log(z)=\int_{[0,z]}\frac{dz}{z}$. Hence we may write \eqref{eq:
scalar Riccati} as
\[
\dot G= A(G-\lambda_+)(G-\lambda_-),\quad G(0,u)=u,
\]
where $\lambda_{\pm}=\frac{-B\pm\sqrt{B^2+4AC}}{2A}$, and it follows
that
\[
G(t,u)=\frac{\lambda_+(u-\lambda_-)-\lambda_-(u-\lambda_+)e^{\lambda
t}}{(u-\lambda_-)-(u-\lambda_+)e^{\lambda t}},\] which can be seen
to be equivalent to \eqref{sol Ric: form1}. As
$\lambda_+\neq\lambda_-$, numerator and denominator cannot vanish at
the same time $t$, and certainly not for $t$ near zero. Hence, by
the maximality of $t_+(u)$, \eqref{sol Ric: form1} is the solution
of \eqref{eq: scalar Riccati} for $t\in [0,t_+(u))$. Finally, the
integral \eqref{eq: integral of G} is checked by differentiation.

\ref{item2: lemsolricc}: We show along the lines of the proof of
Theorem \ref{thmaffine}, that for this choice of coefficients global
solutions exist for initial data $u\in\mathbb C_-$ and stay in
$\mathbb C_-$. To this end, write $R(G)=A G^2 + B G - C$, then
\[
\Re(R(G))=A(\Re(G))^2-A(\Im(G))^2+B\Re(G)-\Re(C)\leq
A(\Re(G))^2+B\Re(G)
\]
and since $A,B\in\mathbb R$ we have that $\Re(G(t,u))\leq 0$ for all
times $t\in [0,t_+(u))$, see Corollary~\ref{corlemstochinvN1} below.
Furthermore, we see that $\Re(\overline G R(G))\leq (1+\vert
G\vert^2)(\vert B\vert+\vert C\vert)$, hence $\partial_t \vert
G(t,u)\vert^2\leq 2(1+\vert G(t,u)\vert^2)(\vert B\vert+\vert
C\vert)$. This implies, by Gronwall's inequality
(\cite[(10.5.1.3)]{dieu60}), that $t_+(u)=\infty$. Hence the lemma
is proved.
\end{proof}

\section{Heston Stochastic Volatility Model}\label{secheston}

This affine model, proposed by Heston \cite{hes93}, generalizes the
Black--Scholes model by assuming a stochastic volatility.

Interest rates are assumed to be constant $r(t)\equiv r\ge 0$, and
there is one risky asset (stock) $S=\e^{X_2}$, where $X=(X_1,X_2)$
is the affine process with state space $\mathbb R_+\times \mathbb R$
and dynamics
\begin{align*}
  dX_1&= (k+\kappa X_1)\,dt +\sigma\sqrt{ 2 X_1}\,dW_1 \\
  dX_2&= (r- X_1)\,dt +\sqrt{2 X_1} \left(\rho\,dW_1+\sqrt{1-\rho^2}
  dW_2\right)
\end{align*}
for some constant parameters $k,\sigma\ge 0$, $\kappa\in\mathbb R$,
and some $\rho\in [-1,1]$. In view of Remark~\ref{remglasskim}, we
note that here \[ \Bcal = \left(\begin{array}{cc} \kappa & 0 \\ -1 &
0\end{array}\right) \] is singular, and hence cannot have strictly
negative eigenvalues.

The implied risk-neutral stock dynamics read
\[ dS = S  r\, dt + S\sqrt{2 X_1}\,d\Wcal\]
for the Brownian motion $\Wcal=\rho  W_1+\sqrt{1-\rho^2}\, W_2$. We
see that $\sqrt{2 X_1}$ is the stochastic volatility of the price
process $S$. They have possibly non-zero covariation
\[ d\langle S,X_1\rangle = 2\rho\sigma S X_1\,dt.\]

The corresponding system of Riccati
equations~\eqref{eqriccfulladmiss} is equivalent to
\begin{equation}\label{ricceqheston}
\begin{aligned}
  \phi(t,u)&=k\int_0^t\psi_1(s,u)\,ds+r u_2 t\\
  \partial_t\psi_1(t,u)&=\sigma^2\psi_1^2(t,u)+(2\rho\sigma u_2+\kappa)\psi_1(t,u)
   +u_2^2  -u_2\\
  \psi_1(0,u)&=u_1\\
  \psi_2(t,u)&=u_2,
\end{aligned}
\end{equation}
which, in view of Lemma \ref{lemsolricc} \ref{item2: lemsolricc}
admits an explicit global solution if $u_1\in\mathbb C_-$ and $0\le
\Re u_2 \le 1$. In particular, for $u_1=0$ and by setting
$\lambda=\sqrt{(2\rho\sigma u_2+\kappa)^2+4\sigma^2(u_2-u_2^2)}$,
the solution can be given explicitly as
\begin{equation}\label{expexpreheston}
\begin{aligned}
\phi(t,u)&=  \frac{k}{\sigma^2}\log\left(\frac{2\lambda
e^{\frac{\lambda-(2\rho\sigma u_2+\kappa)}{2}t}}{\lambda (e^{\lambda
t}+1)-(2\rho\sigma u_2+\kappa)(e^{\lambda t}-1)}\right)+ru_2t
\\
\psi_1(t,u)&= -\frac{2(u_2-u_2^2)(e^{\lambda
t}-1)}{\lambda(e^{\lambda t}+1)-(2\rho\sigma u_2+\kappa)(e^{\lambda
t}-1)}\\ \psi_2(t,u)&= u_2.
\end{aligned}
\end{equation}
Furthermore, for $u=(0,1)$, we obtain
\[ \phi(t,0,1)=rt,\quad\psi(t,0,1) =(0,1)^\top.\]
Theorem~\ref{thmextanaXXX} thus implies that $S(T)$ has finite first
moment, for any $T\in\R_+$, and
\[
\Ex[\e^{-rT}S(T)\mid\Fcal_t]=\e^{-rT}\Ex[\e^{X_2(T)}\mid\Fcal_t]=\e^{-rT}\e^{r(T-t)+X_2(t)}=\e^{-rt}S(t),\]
for $t\le T$, which is just the martingale property of $S$.

We now want to compute the price
\[ \pi(t)=\e^{-r(T-t)} \Ex\left[ (S(T)-K)^+\mid\Fcal_t\right] \]
of a European call option on $S(T)$ with maturity $T$ and strike
price $K$. Fix some ${p}>1$ small enough with $(0,{p})\in
\Dcal_\R(T)$. Formula~\eqref{exoptheston} combined with
\eqref{exopthestonXX} then yields
\begin{multline}\label{eqhestform1}
  \pi(t)=\frac{1}{2\pi}\e^{-r(T-t)} \\
  \times\int_\mathbb R \e^{\phi(T-t,0,{p}+\im y)+\psi_1(T-t,0,{p}+\im y)  X_1(t)+({p}+\im
y)X_2(t)}\,\frac{K^{1-{p}-\im y}}{({p}+\im y)({p}+\im y-1)}\, dy.
\end{multline}
Alternatively, we may fix any $0<{p}<1$ and then, combining
\eqref{eqrelexamf} with \eqref{exopthestonXX},
\begin{multline}\label{eqhestform2}
  \pi(t)=S(t)+\frac{1}{2\pi}\e^{-r(T-t)}\\
  \times\int_\mathbb R \e^{\phi(T-t,0,{p}+\im y)+\psi_1(T-t,0,{p}+\im y)  X_1(t)+({p}+\im
y)X_2(t)}\,\frac{K^{1-{p}-\im y}}{({p}+\im y)({p}+\im y-1)}\, dy.
\end{multline}

Since we have explicit expressions \eqref{expexpreheston} for
$\phi(T-t,0,{p}+\im y)$ and $\psi_1(T-t,0,{p}+\im y)$, we only need
to compute the integral with respect to $y$ in \eqref{eqhestform1}
or \eqref{eqhestform2} numerically.  We have carried out numeric
experiments for European option prices using MATLAB. Fastest results
were achieved for values $p\approx 0.5$ by using \eqref{eqhestform2}
whereas keeping a constant error level the runtime explodes at
$p\rightarrow 0,\,1$, which is due to the singularities of the
integrand. Also, an evaluation of residua
\[
\frac{\pi(t=0,p=1/2)-\pi(t=0,p=1/2+\varepsilon)}{\pi(t=0,p=1/2)}
\]
for $\epsilon\in [0,1/2)\cup(1/2,1]$ suggests that
\eqref{eqhestform2} is numerically more stable than
\eqref{eqhestform1}.

\begin{figure}[h!]
\caption{Implied volatility surface for the Heston
model}\label{fighestonprices}
\centerline{\includegraphics[width=7cm]{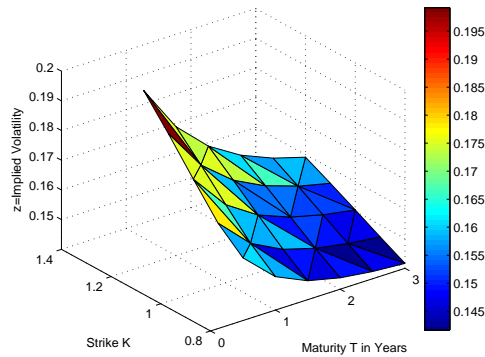}}
\end{figure}

Next, we present implied volatilities obtained by
\eqref{eqhestform2} setting $p=1/2$. As initial data for $X$ and
model parameters, we chose
\[X_1(0) = 0.02,\, X_2(0)=0.00,\, \sigma = 0.1,\, \kappa=-2.0,\,
k=0.02,\, r=0.01,\,\rho=0.5.
\]
Table~\ref{tabhestonprices} shows implied volatilities from call
option prices at $t=0$ for various strikes $K$ and maturities $T$,
computed with \eqref{eqhestform2} for $p=0.5$. These values are in
well accordance with MC simulations (mesh size $T/500$, number of
sample paths $=10000$). The corresponding implied volatility surface
is shown in Figure~\ref{fighestonprices}.

\begin{table}[h!]
\caption{Implied volatilities for the Heston
model}\label{tabhestonprices}
\begin{center}
\begin{tabular}{c*{5}{|l}}
\textit{T-K} &   0.8000  &  0.9000  &  1.0000 &   1.1000 & 1.2000
\\\hline 0.5000  &  0.1611  &  0.1682  &  0.1785  & 0.1892 &
0.1992\\\hline  1.0000  &  0.1513  &  0.1579 &   0.1664 & 0.1751&
0.1835\\\hline  1.5000  &  0.1464  &  0.1524   & 0.1594  & 0.1665&
0.1734\\\hline  2.0000  &  0.1438  &  0.1492  &  0.1551  & 0.1611&
0.1668\\\hline  2.5000  &  0.1424  &  0.1473  &  0.1524  & 0.1574&
0.1623\\\hline  3.0000  &  0.1417  &  0.1460 &   0.1505  & 0.1549&
0.1591\\\hline
\end{tabular}
\end{center}
\end{table}

\begin{remark}
We note that the Heston model is often written in the equivalent
form
\begin{align*}
dv&=\bar{\kappa}(\eta-v)dt+\sigma\sqrt v\, dW_1\\
 dS&=r S dt+ S\sqrt v\,d\mathcal W
\end{align*}
To see the relation of the parameters of this form and the one used
in this section, we simply set $v=2X_1$, and then get
\begin{align*}
dX_1&=(\bar{\kappa}\eta-\bar{\kappa} v )dt+\sigma\sqrt{2X_1}\,
dW_1\quad X_1(0)=X_{10}\\
 \frac{dS}{S}&=r dt+\sqrt{2X_1}\, d\mathcal W,\quad S(0)=e^{X_{2}(0)}
\end{align*}
from which we read off
\[
k=\bar{\kappa}\eta,\quad \kappa=-\bar{\kappa},\quad X_{10}=v_0/2
\]
and all other parameters coincide.

\end{remark}

\section{Affine Transformations and Canonical
Representation}\label{secafftransf}

As above, we let $X$ be affine on the canonical state space $\mathbb
R^m_+\times \mathbb R^n$ with admissible parameters
$a,\alpha_i,b,\beta_i $. Hence, in view of \eqref{sdeaffgen}, for
any $x\in\mathbb R^m_+\times\mathbb R^n$ the process $X=X^x$
satisfies
\begin{equation}\label{sdeX}
\begin{aligned}
 dX &= (b+\Bcal X)\,dt + \rho(X)\,dW\\
 X(0)&=x,
\end{aligned}
\end{equation}
and $\rho(x)\rho(x)^\top=a+\sum_{i\in I} x_i\alpha_i$.

It can easily be checked that for every invertible $d \times
d$-matrix $\Lambda$, the linear transform
\[Y
= \Lambda X  \] satisfies
\begin{equation}\label{Ysdeaff}
dY  = \left(\Lambda b + \Lambda \Bcal \Lambda^{-1} Y  \right)\, dt +
\Lambda \rho\left(\Lambda^{-1}Y \right)\, dW ,\quad Y(0)=\Lambda x .
\end{equation}
 Hence, $Y $ has again an affine drift and diffusion
matrix
\begin{equation} \label{affdiffY}
\Lambda b + \Lambda \Bcal \Lambda^{-1} y  \quad\text{and}\quad
\Lambda\alpha(\Lambda^{-1}y  )\Lambda^\top ,
\end{equation}
respectively.

On the other hand, the affine short rate model~\eqref{r} can be
expressed in terms of $Y(t)$ as
\begin{equation} \label{shrt}
r(t) = c  + \gamma^\top \Lambda^{-1} Y(t) \, .
\end{equation}

This shows that $Y $ and \eqref{shrt} specify an affine short rate
model producing the same short rates, and thus bond prices, as $X $
and \eqref{r}. That is, an invertible linear transformation of the
state process changes the particular form of the stochastic
differential equation \eqref{sdeX}. But it leaves observable
quantities, such as short rates and bond prices invariant.

This motivates the question whether there exists a classification
method ensuring that affine short rate models with the same
observable implications have a unique canonical representation. This
topic has been addressed in
\cite{daisin00,coldufetal06,jos06,chefilkim08}. We now elaborate on
this issue and show that the diffusion matrix $\alpha(x)$ can always
be brought into block-diagonal form by a regular linear transform
$\Lambda$ with $\Lambda(\mathbb R^m_+\times\mathbb R^n)=\mathbb
R^m_+\times\mathbb R^n$.

We denote by
\[{\rm diag}(z_1,\dots,z_m)\] the diagonal matrix with
 diagonal elements $z_1,\dots,z_m$, and we write $I_m$ for the
 $m\times m$-identity matrix.

\begin{lemma}\label{lemcanaffinetrans}
There exists some invertible $d\times d$-matrix $\Lambda$ with
$\Lambda(\mathbb R^m_+\times\mathbb R^n)=\mathbb R^m_+\times\mathbb
R^n$ such that $\Lambda \alpha(\Lambda^{-1}y)\Lambda^\top$ is
block-diagonal of the form
\[\Lambda
\alpha(\Lambda^{-1}y)\Lambda^\top=\left(\begin{array}{cc} {\rm
diag}(y_1,\dots,y_q,0,\dots,0) & 0 \\
 0 & p+\sum_{i\in I} y_i \pi_i\end{array}\right)\]
for some integer $0\le q\le m$ and symmetric positive semi-definite
$n\times n$ matrices
 $p,\pi_1,\dots,\pi_m$. Moreover, $\Lambda b$ and $\Lambda \Bcal \Lambda^{-1}$ meet the
respective admissibility conditions \eqref{eqadmiss} in lieu of $b$
and $\Bcal$.
\end{lemma}

\begin{proof}
From \eqref{eqaffba} we know that $\Lambda\alpha(x)\Lambda^\top$ is
block-diagonal for all $x=\Lambda^{-1} y$ if and only if $\Lambda
a\Lambda^\top$ and $\Lambda\alpha_i\Lambda^\top$ are block-diagonal
for all $i\in I$. By permutation and scaling of the first $m$
coordinate axes (this is a linear bijection from $\mathbb
R^m_+\times\mathbb R^n$ onto itself, which preserves the
admissibility of the transformed $b$ and $\Bcal$), we may assume
that there exists some integer $0\le q\le m$ such that
$\alpha_{1,11}=\cdots=\alpha_{q,qq}=1$ and $\alpha_{i,ii}=0$ for
$q<i\le m$. Hence $a$ and $\alpha_i$ for $q<i\le m$ are already
block-diagonal of the special form
\[ a=\left(\begin{array}{cc} 0 & 0 \\
 0 & a_{JJ}\end{array}\right),\quad \alpha_i=\left(\begin{array}{cc} 0 & 0 \\
 0 & \alpha_{i,JJ}\end{array}\right).\]
For $1\le i\le q$, we may have non-zero off-diagonal elements in the
$i$-th row $\alpha_{i,iJ}$. We thus define the $n\times m$-matrix
${D}=(\delta_1,\dots,\delta_m)$ with $i$-th column
$\delta_{i}=-\alpha_{i,iJ}$ and set
\[ \Lambda=\left(\begin{array}{cc} I_m & 0 \\
 {D} & I_n\end{array}\right).\]
One checks by inspection that ${D}$ is invertible and maps $\mathbb
R^m_+\times\mathbb R^n$ onto $\mathbb R^m_+\times\mathbb R^n$.
Moreover,
\[ {D} \alpha_{i,II} = -\alpha_{i,JI},\quad i\in I.\]
From here we easily verify that
\[ \Lambda \alpha_i=\left(\begin{array}{cc} \alpha_{i,II} & \alpha_{i,IJ} \\
 0 & {D}\alpha_{i,IJ}+\alpha_{i,JJ}\end{array}\right),\]
and thus
\[ \Lambda \alpha_i\Lambda ^\top=\left(\begin{array}{cc} \alpha_{i,II} & 0 \\
 0 & {D}\alpha_{i,IJ}+\alpha_{i,JJ}\end{array}\right).\]
Since $\Lambda a \Lambda^\top =a$, the first assertion is proved.

The admissibility conditions for $\Lambda b$ and $\Lambda \Bcal
\Lambda^{-1}$ can easily be checked as well.
\end{proof}

In view of \eqref{affdiffY}, \eqref{shrt} and
Lemma~\ref{lemcanaffinetrans} we thus obtain the following result.

\begin{theorem}[Canonical Representation]\label{thmcanaff}
Any affine short rate model \eqref{r}, after some modification of
$\gamma$ if necessary, admits an $\mathbb R^m_+\times \mathbb
R^n$-valued affine state process $X$ with block-diagonal diffusion
matrix of the form
\begin{equation}\label{eqcanaffmat}
\alpha(x)=\left(\begin{array}{cc} {\rm
diag}(x_1,\dots,x_q,0,\dots,0) & 0 \\
 0 & a+\sum_{i\in I} x_i \alpha_{i,JJ}\end{array}\right)
\end{equation}
for some integer $0\le q\le m$.
\end{theorem}

\section{Existence and Uniqueness of Affine
Processes}\label{secexiuniaffine}

All we said about the affine process $X$ so far was under the
premise that there exists a unique solution $X=X^x$ of the
stochastic differential equation \eqref{sdeaffgen} on some
appropriate state space $\Xcal\subset\mathbb R^d$. However, if the
diffusion matrix $\rho(x)\rho(x)^\top$ is affine then $\rho(x)$
cannot be Lipschitz continuous in $x$ in general. This raises the
question whether \eqref{sdeaffgen} admits a solution at all.

In this section, we show how $X$ can always be realized as unique
solution of the stochastic differential equation \eqref{sdeaffgen},
which is \eqref{sdeX}, in the canonical affine framework
$\Xcal=\mathbb R^m_+\times\mathbb R^n$ and for particular choices of
$\rho(x)$.

We recall from Theorem~\ref{thmaffchar} that the affine property of
$X$ imposes explicit conditions on $\rho(x)\rho(x)^\top$, but not on
$\rho(x)$ as such. Indeed, for any orthogonal $d\times d$-matrix
$D$, the function $\rho(x)D$ yields the same diffusion matrix,
$\rho(x)D D^\top \rho(x)^\top=\rho(x) \rho(x)^\top$, as $\rho(x)$.

On the other hand, from Theorem~\ref{thmaffine} we know that any
admissible parameters $a,\alpha_i,b,\beta_i$ in \eqref{eqaffba}
uniquely determine the functions
$(\phi(\cdot,u),\psi(\cdot,u)):\mathbb R_+\to \mathbb
C_-\times\mathbb C^m_-\times\im\mathbb R^n$ as solution of the
Riccati equations~\eqref{eqriccfulladmiss}, for all $u\in\mathbb
C^m_-\times\im\mathbb R^n$. These in turn uniquely determine the law
of the process $X$. Indeed, for any $0\le t_1<t_2$ and
$u_1,\,u_2\in\mathbb C^m_-\times\im\mathbb R^n$, we infer by
iteration of \eqref{eqaffdef1}

\begin{align*}
 \Ex\left[ \e^{u_1^\top  X(t_1) + u_2^\top  X(t_2)}\right] &=
 \Ex\left[ \e^{u_1^\top  X(t_1)} \Ex\left[ \e^{u_2^\top
X(t_2)}\mid\Fcal_{t_1}\right]\right] \\
 &=
\Ex\left[ \e^{u_1^\top
X(t_1)}\e^{\phi(t_2-t_1,u_2)+\psi(t_2-t_1,u_2)^\top X(t_1)}\right]\\
&=\e^{\phi(t_2-t_1,u_2)+\phi(t_1,u_1+\psi(t_2-t_1,u_2))+\psi(t_1,u_1+\psi(t_2-t_1,u_2))^\top
x}.
\end{align*}
Hence the joint distribution of $(X(t_1),X(t_2))$ is uniquely
determined by the functions $\phi$ and $\psi$. By further iteration
of this argument, we conclude that every finite dimensional
distribution, and thus the law, of $X$ is uniquely determined by the
parameters $a,\alpha_i,b,\beta_i$.

We conclude that the law of an affine process $X$, while uniquely
determined by its characteristics \eqref{eqaffba}, can be realized
by infinitely many variants of the stochastic differential equation
\eqref{sdeX} by replacing $\rho(x)$ by $ \rho(x) D$, for any
orthogonal $d\times d$-matrix $D$. We now propose a canonical choice
of $\rho(x)$ as follows:

\begin{itemize}

\item In view of \eqref{Ysdeaff} and Lemma~\ref{lemcanaffinetrans}, every affine process
$X$ on $\mathbb R^m_+\times\mathbb R^n$ can be written as
$X=\Lambda^{-1} Y$ for some invertible $d\times d$-matrix $\Lambda$
and some affine process $Y$ on $\mathbb R^m_+\times\mathbb R^n$ with
block-diagonal diffusion matrix. It is thus enough to consider such
$\rho(x)$ where $\rho(x)\rho(x)^\top$ is of the form
\eqref{eqcanaffmat}. Obviously, $\rho(x)\equiv \rho(x_I)$ is a
function of $x_I$ only.

\item Set $\rho_{IJ}(x)\equiv 0$, $\rho_{JI}(x)\equiv 0$, and
\[\rho_{II}(x_I)= {\rm diag}(\sqrt{x_1},\dots,\sqrt{x_q},0,\dots,0) .\]
Chose for $\rho_{JJ}(x_I)$ any measurable $n\times n$-matrix-valued
function satisfying
\begin{equation}\label{eqrhoJJ}
\rho_{JJ}(x_I) \rho_{JJ}(x_I)^\top=a+\sum_{i\in I} x_i\alpha_{i,JJ}
.
\end{equation}
In practice, one would determine $\rho_{JJ}(x_I)$ via Cholesky
factorization, see e.g. \cite[Theorem~2.2.5]{neu01}. If
$a+\sum_{i\in I} x_i\alpha_{i,JJ}$ is strictly positive definite,
then $\rho_{JJ}(x_I)$ turns out to be the unique lower triangular
matrix with strictly positive diagonal elements and satisfying
\eqref{eqrhoJJ}. If $a+\sum_{i\in I} x_i\alpha_{i,JJ}$ is merely
positive semi-definite, then the algorithm becomes more involved. In
any case, $\rho_{JJ}(x_I)$ will depend measurably on $x_I$.
\item The stochastic differential equation~\eqref{sdeX} now reads
\begin{equation}\label{sdeaffcanX}
 \begin{aligned}
  dX_I&= (b_I+\Bcal_{II} X_I)\,dt + \rho_{II}(X_I)\,dW_I\\
dX_J&= (b_J+\Bcal_{JI} X_I + \Bcal_{JJ} X_J)\,dt +
\rho_{JJ}(X_I)\,dW_J\\ X(0)&=x
\end{aligned}
\end{equation}
Lemma~\ref{lemexiuniyamwat} below asserts the existence and
uniqueness of an $\mathbb R^m_+\times\mathbb R^n$-valued solution
$X=X^x$, for any $x\in\mathbb R^m_+\times\mathbb R^n$.
\end{itemize}

We thus have shown:
\begin{theorem}
Let $a,\alpha_i,b,\beta_i$ be admissible parameters. Then there
exists a measurable function $\rho:\mathbb R^m_+\times\mathbb R^n\to
\mathbb R^{d\times d}$ with $\rho(x)\rho(x)^\top=a+\sum_{i\in I}
x_i\alpha_i$, and such that, for any $x\in\mathbb R^m_+\times\mathbb
R^n$, there exists a unique $\mathbb R^m_+\times\mathbb R^n$-valued
solution $X=X^x$ of \eqref{sdeX}.

Moreover, the law of $X$ is uniquely determined by
$a,\alpha_i,b,\beta_i$, and does not depend on the particular choice
of $\rho$.
\end{theorem}

The proof of the following lemma uses the concept of a weak
solution. The interested reader will find detailed background in
e.g.\ \cite[Section 5.3]{kar/shr:91}.

\begin{lemma}\label{lemexiuniyamwat}
For any $x\in\mathbb R^m_+\times\mathbb R^n$, there exists a unique
$\mathbb R^m_+\times\mathbb R^n$-valued solution $X=X^x$ of
\eqref{sdeaffcanX}.
\end{lemma}

\begin{proof}
First, we extend $\rho$ continuously to $\mathbb R^d$ by setting
$\rho(x)=\rho(x_1^+,\dots,x_m^+)$, where we denote
$x_i^+=\max(0,x_i)$.

Now observe that $X_I$ solves the autonomous equation
\begin{equation}\label{sdeXI}
dX_I = (b_I+\Bcal_{II} X_I)\,dt + \rho_{II}(X_I)\,dW_I,\quad
X_I(0)=x_I.
\end{equation}
Obviously, there exists a finite constant $K$ such that the linear
growth condition
\[ \|b_I+\Bcal_{II} x_I\|^2+\|\rho(x_I)\|^2 \le K(1+\|x_I\|^2) \]
is satisfied for all $x\in\mathbb R^m$. By \cite[Theorems 2.3 and
2.4]{ikewat81} there exists a weak solution\footnote{A weak solution
consists of a filtered probability space
$(\Omega,\Fcal,(\Fcal_t),\Pa)$ carrying a continuous adapted process
$X_I$ and a Brownian motion $W_I$ such that \eqref{sdeXI} is
satisfied. The crux of a weak solution is that $X_I$ is not
necessarily adapted to the filtration generated by the Brownian
motion $W_I$. See \cite[Definition 1]{yamwat71} or \cite[Definition
5.3.1]{kar/shr:91}.} of \eqref{sdeXI}. On the other hand,
\eqref{sdeXI} is exactly of the form as assumed in \cite[Theorem
1]{yamwat71}, which implies that pathwise
uniqueness\footnote{Pathwise uniqueness holds if, for any two weak
solutions $(X_I,W_I)$ and $(X_I',W_I)$ of \eqref{sdeXI} defined on
the the same probability space $(\Omega,\Fcal,\Pa)$ with common
Brownian motion $W_I$ and with common initial value
$X_I(0)=X_I'(0)$, the two processes are indistinguishable:
$\Pa[X_I(t)=X_I'(t)$ for all $t\ge 0]=1$. See \cite[Definition
2]{yamwat71} or \cite[Section 5.3]{kar/shr:91}.} holds for
\eqref{sdeXI}. The Yamada--Watanabe Theorem, see \cite[Corollary
3]{yamwat71} or \cite[Corollary 5.3.23]{kar/shr:91}, thus implies
that there exists a unique solution $X_I=X_I^{x_I}$ of
\eqref{sdeXI}, for all $x_I\in\mathbb R^m$.

Given $X_I^{x_I}$, it is then easily seen that
\begin{align*}
X_J(t)&=\e^{\Bcal_{JJ} t} \left( x_J + \int_0^t \e^{-\Bcal_{JJ}
s}(b_J+ \Bcal_{JI} X_I(s))\,ds \right.\\&+\left. \int_0^t
\e^{-\Bcal_{JJ} s}\rho_{JJ}(X_I(s))\,dW_J(s) \right)
\end{align*} is the unique solution to the second equation in
\eqref{sdeaffcanX}.

Admissibility of the parameters $b$ and $\beta_i$ and the stochastic
invariance Lemma~\ref{lemstochinvN} eventually imply that
$X_I=X_I^{x_I}$ is $\mathbb R^m_+$-valued for all $x_I\in\mathbb
R^m_+$. Whence the lemma is proved.
\end{proof}

\begin{appendix}

\section{On the Regularity of Characteristic
Functions}\label{secregcharfun}

This auxiliary section provides some analytic regularity results for
characteristic functions, which are of independent interest. These
results enter the main text only via the proof of
Theorem~\ref{thmextanaXXX}. This section may thus be skipped at the
first reading.

Let $\nu$ be a bounded measure on $\mathbb R^d$, and denote by
\[ G(z)=\int_{\mathbb R^d} \e^{ z^\top x}\, \nu(dx) \]
its characteristic function\footnote{This is a slight abuse of
terminology, since the characteristic function $g(y)=G(\im y)$ of
$\nu$ is usually defined on real arguments $y\in\mathbb R^d$.
However, it facilitates the subsequent notation.} for $z\in
\im\mathbb R^d$. Note that $G(z)$ is actually well defined for $z\in
\Scal(V)$ where
\[ V =\left\{ y\in\mathbb R^d\mid \int_{\mathbb R^d} \e^{  y^\top x}\,
\nu(dx)<\infty\right\}.\]

We first investigate the interplay between the (marginal) moments of
$\nu$ and the corresponding (partial) regularity of $G$.

\begin{lemma}\label{lemRC1}
Denote $g(y)=G(iy)$ for $y\in\mathbb R^d$, and let $k\in\mathbb N$
and $1\le i\le d$.

If $\partial_{y_i}^{2k} g(0)$
  exists then
  \[ \int_{\mathbb R^d} |x_i|^{2k}\,\nu(dx)<\infty.\]
On the other hand, if $\int_{\mathbb R^d} \| x\|^k\,\nu(dx)<\infty$
then
  $g\in C^k$ and
  \[ \partial_{y_{i_1}}\cdots\partial_{y_{i_l}}
  g(y)=\im^l\int_{\mathbb R^d} x_{i_1}\cdots x_{i_l}\,\e^{\im y^\top
  x}\,\nu(dx)\]
  for all $y\in\mathbb R^d$, $1\le i_1,\dots,i_l\le d$ and $1\le l\le k$.
\end{lemma}

\begin{proof}
As usual, let $e_i$ denote the $i$th standard basis vector in
$\mathbb R^d$. Observe that $s\mapsto g(s e_i)$ is the
characteristic function of the image measure of $\nu$ on $\mathbb R$
by the mapping $x\mapsto x_i$. Since $\partial_s^{2k}
g(se_i)|_{s=0}=\partial_{y_i}^{2k} g(0)$, the assertion follows from
the one-dimensional case, see \cite[Theorem 2.3.1]{luk70}.

The second part of the lemma follows by differentiating under the
integral sign, which is allowed by dominated convergence.
\end{proof}

\begin{lemma}\label{lemRC2X}
The set $V$ is convex. Moreover, if $U\subset V $ is an open set in
  $\mathbb R^d$, then $G$ is analytic on the open strip $\Scal(U)$ in $\mathbb C^d$.
\end{lemma}

\begin{proof}
Since $G:\mathbb R^d\to [0,\infty]$ is a convex function, its domain
$V=\{y\in\mathbb R^d\mid G(y)<\infty\}$ is convex, and so is every
level set $V_l=\{y\in\mathbb R^d\mid G(y)\le l\}$ for $l\ge 0$.

Now let $U\subset V $ be an open set in $\mathbb R^d$. Since any
convex function on $\mathbb R^d$ is continuous on the open interior
of its domain, see \cite[Theorem 10.1]{rock}, we infer that $G$ is
continuous on $U$. We may thus assume that $U_l=\{y\in\mathbb
R^d\mid G(y)<l\}\cap U\subset V_l$ is open in $\mathbb R^d$ and
non-empty for $l>0$ large enough.

Let $z\in \Scal(U_l)$ and $(z_n)$ be a sequence in $\Scal(U_l)$ with
$z_n\to z$. For $n$ large enough, there exists some $p>1$ such that
$p z_n\in \Scal(U_l)$. This implies $p\mathbb Re z_n\in V_l$ and
hence
\[ \int_{\mathbb R^d} \left| \e^{z_n^\top x}\right|^p\,\nu(dx)\le l .\]
Hence the class of functions $\{\e^{z_n^\top x}\mid n\in\mathbb N\}$
is uniformly integrable with respect to $\nu$, see
\cite[13.3]{will91}. Since $\e^{z_n^\top x}\to \e^{z^\top x}$ for
all $x$, we conclude by Lebesgue's convergence theorem that
\[ |G(z_n)-G(z)|\le\int_{\mathbb R^d} \left| \e^{z_n^\top x}-\e^{z^\top
x}\right| \,\nu(dx)\to 0.\] Hence $G$ is continuous on $\Scal(U_l)$.

It thus follows from the Cauchy formula, see \cite[Section
IX.9]{dieu60}, that $G$ is analytic on $\Scal(U_l)$ if and only if,
for every $z\in\Scal(U_l)$ and $1\le i\le d$, the function
$\zeta\mapsto G(z+\zeta e_i)$ is analytic on $\{\zeta\in\mathbb
C\mid z+\zeta e_i\in\Scal(U_l)\}$. Here, as usual, we denote $e_i$
the $i$th standard basis vector in $\mathbb R^d$.

We thus let $z\in\Scal(U_l)$ and $1\le i\le d$. Then there exists
some $\epsilon_-<0<\epsilon_+$ such that $  z+ \zeta e_i\in
\Scal(U_l)$ for all $\zeta\in\Scal( [\epsilon_-,\epsilon_+])$. In
particular, $|\e^{(z+\epsilon_- e_i)^\top x}|\,\nu(dx)$ and
$|\e^{(z+\epsilon_+ e_i)^\top x}|\,\nu(dx)$ are bounded measures on
$\mathbb R^d$. By dominated convergence, it follows that the two
summands
\begin{align*}
G(z+\zeta e_i)&=\int_{\{ x_i<0\}} \e^{(\zeta-\epsilon_-)
x_i}\,\e^{(z+\epsilon_- e_i)^\top x}\,\nu(dx) \\&+\int_{\{ x_i\ge
0\}} \e^{(\zeta-\epsilon_+)x_i}\,\e^{(z+\epsilon_+ e_i)^\top
x}\,\nu(dx),
\end{align*}
are complex differentiable, and thus $G$ is analytic, in
$\zeta\in\Scal((\epsilon_-,\epsilon_+))$. Whence $G$ is analytic on
$\Scal(U_l)$. Since $\Scal(U)=\cup_{l>0}\Scal(U_l)$, the lemma
follows.
\end{proof}

In general, $V$ does not have an open interior in $\mathbb R^d$. The
next lemma provides sufficient conditions for the existence of an
open set $U\subset V$ in $\mathbb R^d$.

\begin{lemma}\label{lemRC3X}
Let $U'$ be an open neighborhood of $0$ in $\mathbb C^d$ and $h$ an
analytic function on $U'$. Suppose that $U=U'\cap\mathbb R^d$ is
star-shaped around $0$ and $G(z)=h(z)$ for all $z\in
U'\cap\im\mathbb R^d$. Then $U \subset V$ and $G=h$ on $U'\cap
\Scal(U)$.
\end{lemma}

\begin{proof}

We first suppose that $U'=P_\rho$ for the open polydisc
\[ P_\rho =\left\{ z\in\mathbb C^d\mid |z_i |<\rho_i,\;1\le
i\le d\right\},\] for some $\rho=(\rho_1,\dots,\rho_d)\in\mathbb
R^d_{++}$. Note the symmetry $\im P_\rho=P_\rho$.

As in Lemma~\ref{lemRC1}, we denote $g(y)=G(iy)$ for $y\in\mathbb
R^d$. By assumption, $g(y)=h(\im y)$ for all $y\in P_\rho\cap\mathbb
R^d$. Hence $g$ is analytic on $P_{\rho}\cap\mathbb R^d$, and the
Cauchy formula, \cite[Section IX.9]{dieu60}, yields
\[ g(y)=\sum_{i_1,\dots,i_d\in\mathbb N_0} c_{i_1,\dots,i_d}
 y_1^{i_1}\cdots y_d^{i_d}\quad \text{for $y\in P_\rho\cap\mathbb R^d$}\] where
$\sum_{i_1,\dots,i_d\in\mathbb N_0} c_{i_1,\dots,i_d}
 z_1^{i_1}\cdots z_d^{i_d}=h(\im z)$ for all $z\in P_{\rho}$.
This power series is absolutely convergent on $P_{\rho}$, that is,
\[ \sum_{i_1,\dots,i_d\in\mathbb N_0} |c_{i_1,\dots,i_d}|\,
 |z_1^{i_1}\cdots z_d^{i_d}|<\infty\quad \text{for all $z\in
 P_{\rho}$.}\]

From the first part of Lemma~\ref{lemRC1}, we infer that $\nu$
possesses all moments, that is, $\int_{\mathbb
R^d}\|x\|^k\,\nu(dx)<\infty$ for all $k\in\mathbb N$. From the
second part of Lemma~\ref{lemRC1} thus
\[ c_{i_1,\dots,i_d}=\frac{\im^{i_1+\cdots+i_d}}{i_1!\cdots
i_d!}\int_{\mathbb R^d} x_1^{i_1}\cdots x_d^{i_d}\,\nu(dx).\]

From the inequality $|x_i|^{2k-1}\le (x_i^{2k}+x^{2k-2}_i)/2$, for
$k\in\mathbb N$, and the above properties, we infer that for all
$z\in P_{\rho}$,
\[\int_{\mathbb R^d} \e^{\sum_{i=1}^d
|z_i|\,|x_i|}\,\nu(dx)=\sum_{i_1,\dots,i_d\in\mathbb N_0} \frac{
|z_1^{i_1}\cdots z_d^{i_d}|}{i_1!\cdots i_d!}\int_{\mathbb R^d}
|x_1^{i_1}\cdots x_d^{i_d}|\,\nu(dx)<\infty
\] Hence
$P_{\rho}\cap\mathbb R^d\subset V$, and Lemma~\ref{lemRC2X} implies
that $G$ is analytic on $\Scal(P_{\rho}\cap\mathbb R^d)$. Since the
power series for $G$ and $h$ coincide on $P_\rho\cap\im\mathbb R^d$,
we conclude that $G=h$ on $P_\rho$, and the lemma is proved for
$U'=P_\rho$.

Now let $U'$ be an open neighborhood of $0$ in $\mathbb C^d$. Then
there exists some open polydisc $P_\rho\subset U'$ with
$\rho\in\mathbb R^d_{++}$. By the preceding case, we have
$P_\rho\cap\mathbb R^d\subset V$ and $G=h$ on $P_\rho$. In view of
Lemma~\ref{lemRC2X} it thus remains to show that $U=U'\cap\mathbb
R^d\subset V$.

To this end, let $a\in U$. Since $U$ is star-shaped around $0$ in
$\R^d$, there exists some $s_1>1$ such that $sa\in U$ for all $s\in
[0,s_1]$ and $h(sa)$ is analytic in $s\in (0,s_1)$. On the other
hand, there exists some $0<s_0<s_1$ such that $sa\in P_\rho\cap\R^d$
for all $s\in [0,s_0]$, and $G(sa)=h(sa)$ for $s\in (0,s_0)$. This
implies
\[  \int_{\{ a^\top x\ge 0\}} \e^{s a^\top x}\,\nu(dx)=h(sa)-\int_{\{ a^\top x< 0\}}
\e^{s a^\top x}\,\nu(dx)\] for $s\in (0,s_0)$. By
Lemma~\ref{lemRC2X}, the right hand side is an analytic function in
$s\in (0,s_1)$. We conclude by Lemma~\ref{lemmgfRp} below, for $\mu$
defined as the image measure of $\nu$ on $\R_+$ by the mapping
$x\mapsto a^\top x$, that $a\in V$. Hence the lemma is proved.
\end{proof}

\begin{lemma}\label{lemmgfRp}
Let $\mu$ be a bounded measure on $\R_+$, and $h$ an analytic
function on $(0,s_1)$, such that
\begin{equation}\label{lemmgfRpeq1}
  \int_{\R_+} \e^{sx}\,\mu(dx)=h(s)
\end{equation}
for all $s\in (0,s_0)$, for some numbers $0<s_0<s_1$. Then
\eqref{lemmgfRpeq1} also holds for $s\in (0,s_1)$.
\end{lemma}

\begin{proof}
Denote $f(s)=\int_{\R_+} \e^{sx}\,\mu(dx)$ and define
$s_\infty=\sup\left\{s>0\mid f(s)<\infty\right\}\ge s_0$, such that
\begin{equation}\label{lemmgfRpeq2}
 f(s)=+\infty\quad\text{for $s\ge s_\infty$.}
\end{equation}

We assume, by contradiction, that $s_\infty<s_1$. Then there exists
some $s_\ast\in (0,s_\infty)$ and $\epsilon>0$ such that
$s_\ast<s_\infty<s_\ast+\epsilon$ and such that $h$ can be developed
in an absolutely convergent power series
\[ h(s) = \sum_{k\ge 0} \frac{c_k}{k!} (s-s_\ast)^k  \quad \text{for
$s\in (  s_\ast-\epsilon,  s_\ast+\epsilon)$.}\] In view of
Lemma~\ref{lemRC2X}, $f$ is analytic, and thus $f=h$, on
$(0,s_\infty)$. Hence we obtain, by dominated convergence,
\[ c_k=\frac{d^k}{ds^k}h(s)|_{s=s_\ast} =
\frac{d^k}{ds^k}f(s)|_{s=s_\ast}=\int_{\R_+} x^k\e^{s_\ast
x}\,\mu(dx)\ge 0.\] By monotone convergence, we conclude
\[ h(s) = \sum_{k\ge 0} \int_{\R_+} \frac{x^k}{k!}(s-s_\ast)^k \e^{s_\ast
x}\,\mu(dx)= \int_{\R_+} \sum_{k\ge 0}\frac{x^k}{k!}(s-s_\ast)^k
\e^{s_\ast x}\,\mu(dx) =\int_{\R_+} \e^{sx}\,\mu(dx)\] for all $s\in
(s_\ast,s_\ast+\epsilon)$. But this contradicts \eqref{lemmgfRpeq2}.
Whence $s_\infty\ge s_1$, and the lemma is proved.
\end{proof}

\section{Invariance and Comparison Results for Differential
Equations}\label{sec: comp}

In this section we deliver invariance and comparison results for
stochastic and ordinary differential equations, which are used in
the proofs of the main Theorems~\ref{thmaffine}, \ref{thmextanaXXX}
and \ref{thmaffdisc1} and Lemma~\ref{lemexiuniyamwat} above.

We start with an invariance result for the stochastic differential
equation~\eqref{sdeaffgen}.

\begin{lemma}\label{lemstochinvN}
Suppose $b$ and $\rho$ in \eqref{sdeaffgen} admit a continuous and
measurable extension to $\mathbb R^d$, respectively, and such that
$a$ is continuous on $\mathbb R^d$. Let $u\in\mathbb
R^d\setminus\{0\}$ and define the half space
\[ H=\{ x\in\mathbb R^d\mid u^\top x\ge 0\},\]
its interior $H^0=\{ x\in\mathbb R^d\mid u^\top x> 0\}$, and its
boundary $\partial H=\{ x\in H\mid u^\top x=0\}$.

\begin{enumerate}
  \item\label{lemstochinvNa} Fix $x\in\partial H$ and let $X=X^x$ be a solution of
\eqref{sdeaffgen}. If $X(t)\in H$ for all $t\ge 0$, then necessarily
\begin{align}
  u^\top a(x)\,u &= 0 \label{eqdiffpar}\\
  u^\top b(x) &\ge 0.\label{eqdriftin}
\end{align}

\item\label{lemstochinvNb} Conversely, if \eqref{eqdiffpar} and \eqref{eqdriftin} hold
for all $x\in \mathbb R^d\setminus H^0$, then any solution $X$ of
\eqref{sdeaffgen} with $X(0)\in H$ satisfies $X(t)\in H$ for all
$t\ge 0$.
\end{enumerate}
\end{lemma}

Intuitively speaking, \eqref{eqdiffpar} means that the diffusion
must be ``parallel to the boundary'', and \eqref{eqdriftin} says
that the drift must be ``inward pointing'' at the boundary of $H$.

\begin{proof}
Fix $x\in\partial H$ and let $X=X^x$ be a solution of
\eqref{sdeaffgen}. Hence
\[ u^\top X(t) = \int_0^t u^\top b(X(s))\,ds +\int_0^t
u^\top\rho(X(s))\,dW(s).\] Since $a$ and $b$ are continuous, there
exists a stopping time $\tau_1>0$ and a finite constant $K$ such
that
\[ |u^\top b(X(t\wedge\tau_1))|\le K \]
and \[ \|u^\top\rho(X(t\wedge\tau_1))\|^2=u^\top
a(X(t\wedge\tau_1))\,u \le K \] for all $t\ge 0$. In particular, the
stochastic integral part of $u^\top X(t\wedge\tau_1)$ is a
martingale. Hence
\[ \Ex\left[ u^\top X(t\wedge\tau_1) \right]=\Ex\left[ \int_0^{t\wedge\tau_1} u^\top
b(X(s))\,ds \right],\quad t\ge 0.\]

We now argue by contradiction, and assume first that $u^\top
b(x)<0$. By continuity of $b$ and $X(t)$, there exists some
$\epsilon>0$ and a stopping time $\tau_2>0$ such that $u^\top
b(X(t))\le -\epsilon$ for all $t\le \tau_2$. In view of the above
this implies
\[ \Ex\left[ u^\top X(\tau_2\wedge\tau_1) \right]<0.\]
This contradicts $X(t)\in H$ for all $t\ge 0$, whence
\eqref{eqdriftin} holds.

As for \eqref{eqdiffpar}, let $C>0$ be a finite constant and define
the stochastic exponential $Z_t =\Ecal (-C\int_0^tu^\top\rho(X
)\,dW)$. Then $Z$ is a strictly positive local martingale.
Integration by parts yields
\[ u^\top X(t) Z(t) = \int_0^t Z(s)\left(u^\top b(X(s))-C\,u^\top
a(X(s))\,u\right)ds + M(t)\] where $M$ is a local martingale. Hence
there exists a stopping time $\tau_3>0$ such that for all $t\ge 0$,
\[ \Ex\left[ u^\top X(t\wedge\tau_3) Z(t\wedge\tau_3) \right]
=\Ex\left[ \int_0^{t\wedge\tau_3} Z(s)\left(u^\top b(X(s))-C\,u^\top
a(X(s))\,u\right)ds \right].\] Now assume that $u^\top a(x)\,u 0$.
By continuity of $a$ and $X(t)$, there exists some $\epsilon>0$ and
a stopping time $\tau_4>0$ such that $u^\top a(X(t))\,u\ge\epsilon$
for all $t\le \tau_4$. For $C>K/\epsilon$, this implies
\[ \Ex\left[ u^\top X(\tau_4\wedge\tau_3\wedge\tau_1) Z(\tau_4\wedge\tau_3\wedge\tau_1)
\right]<0.\] This contradicts $X(t)\in H$ for all $t\ge 0$. Hence
\eqref{eqdiffpar} holds, and part~{\ref{lemstochinvNa}} is proved.

As for part~\ref{lemstochinvNb}, suppose \eqref{eqdiffpar} and
\eqref{eqdriftin} hold for all $x\in \mathbb R^d\setminus H^0$, and
let $X$ be a solution of \eqref{sdeaffgen} with $X(0)\in H$. For
$\delta,\epsilon>0$ define the stopping time
\[ \tau_{\delta,\epsilon}=\inf\left\{ t\mid \text{$u^\top
X(t)\le-\epsilon$ and $u^\top X(s)<0$ for all $s\in
[t-\delta,t]$}\right\}.\] Then on
$\{\tau_{\delta,\epsilon}<\infty\}$ we have $u^\top\rho(X(s))=0$ for
$\tau_{\delta,\epsilon}-\delta\le s\le \tau_{\delta,\epsilon}$ and
thus
\[ 0>u^\top X(\tau_{\delta,\epsilon})-u^\top
X(\tau_{\delta,\epsilon}-\delta)=\int_{\tau_{\delta,\epsilon}-\delta}^{\tau_{\delta,\epsilon}}
u^\top b(X(s))\,ds\ge 0,\] a contradiction. Hence
$\tau_{\delta,\epsilon}=\infty$. Since $\delta,\epsilon>0$ were
arbitrary, we conclude that $u^\top X(t)\ge 0$ for all $t\ge 0$, as
desired. Whence the lemma is proved.
\end{proof}

It is straightforward to extend Lemma~\ref{lemstochinvN} towards a
polyhedral convex set $\cap_{i=1}^k H_i$ with half-spaces $H_i=\{
x\in\mathbb R^d\mid u_i^\top x\ge 0\}$, for some elements
$u_1,\dots,u_k\in\mathbb R^d\setminus\{0\}$ and some $k\in\mathbb
N$. This holds in particular for the canonical state space $\mathbb
R^m_+\times\mathbb R^n$. Moreover, Lemma~\ref{lemstochinvN} includes
time-inhomogeneous\footnote{Time-inhomogeneous differential
equations can be made homogeneous by enlarging the state space.}
ordinary differential equations as special case. The proofs of the
following two corollaries are left to the reader.

\begin{corollary}\label{corlemstochinvN1}
  Let $H_i=\{
x\in\mathbb R^d\mid x_i\ge 0\}$ denote the $i$-th canonical half
space in $\R^d$, for $i=1,\dots,m$. Let $b:\R_+\times \R^d\to\R^d$
be a continuous map satisfying, for all $t\ge 0$,
\begin{align*}
    b(t,x)&=b(t,x_1^+,\dots,x_m^+,x_{m+1},\dots,x_d)\quad\text{for all $x\in \R^d$, and}\\
    b_i(t,x)&\ge 0\quad \text{for all $x\in \partial H_i$, $i=1,\dots,m$.}
\end{align*}
Then any solution $f$ of
\[ \partial_t f(t)=b(t,f(t)) \]
with $f(0)\in \R^m_+ \times \R^n$ satisfies $f(t)\in
\R^m_+\times\R^n$ for all $t\ge 0$.
\end{corollary}

\begin{corollary}\label{corlincomp}
Let $B(t)$ and $C(t)$ be continuous $\R^{m\times m}$- and
$\R^m_+$-valued parameters, respectively, such that $B_{ij}(t)\ge 0$
whenever $i\neq j$. Then the solution $f$ of the linear differential
equation in $\R^m$
\[\partial_t f(t)=B(t)\,f(t)+C(t)\]
with $f(0)\in \mathbb R^m_+$ satisfies $f(t)\in\R^m_+$ for all
$t\geq 0$.
\end{corollary}

Here and subsequently, we let $\succeq$ denote the partial order on
$\R^m$ induced by the cone $\R^m_+$. That is, $x\succeq y$ if
$x-y\in\R^m_+$. Then Corollary~\ref{corlincomp} may be rephrased,
for $C(t)\equiv 0$, by saying that the operator $\e^{\int_0^t
B(s)\,ds}$ is $\succeq$-order preserving, i.e.\ $\e^{\int_0^t
B(s)\,ds}\R^m_+\subseteq\R^m_+$.

Next, we consider time-inhomogeneous Riccati equations in $\R^m$ of
the special form
\begin{equation}\label{ricccomp}
  \partial_t f_i(t)=A_i  f_i(t)^2+B_i^\top f(t)+C_i(t) ,\quad i=1,\dots,m,
\end{equation}
for some parameters $A,B,C(t)$ satisfying the following
admissibility conditions
\begin{equation}\label{adminparNN}
\begin{gathered}
  A =(A_1 ,\dots,A_m )\in\R^m,\\
  B_{i,j} \ge 0\quad\text{for $1\le i\neq j\le m$,}\\
  \text{$C(t)=(C_1(t),\dots,C_m(t))$ continuous $\R^m$-valued.}
\end{gathered}
\end{equation}

The following lemma provides a comparison result for
\eqref{ricccomp}. It shows, in particular, that the solution of
\eqref{ricccomp} is uniformly bounded from below on compacts with
respect to $\succeq$ if $A\succeq 0$.

\begin{lemma}\label{lemcompriccX}
Let $A^{(k)},B ,C^{(k)}$, $k=1,2$, be parameters satisfying the
admissibility conditions~\eqref{adminparNN}, and
\begin{equation}
  A^{(1)}\preceq A^{(2)},\quad C^{(1)}(t)\preceq C^{(2)}(t).
\end{equation}
Let $\tau>0$ and $f^{(k)}:[0,\tau)\to\R^m$ be solutions of
\eqref{adminparNN} with $A$ and $C$ replaced by $A^{(k)}$ and
$C^{(k)}$, respectively, $k=1,2$. If $f^{(1)}(0)\preceq  f^{(2)}(0)$
then $f^{(1)}(t)\preceq f^{(2)}(t)$ for all $t\in [0,\tau)$. If,
moreover, $A^{(1)}=0$ then
\[ \e^{B t} \left(f^{(1)}(0)+\int_0^t \e^{-B s} C^{(1)}(s)\,ds\right)\preceq
f^{(2)}(t)\] for all $t\in [0,\tau)$.
\end{lemma}

\begin{proof}
The function $f=f^{(2)}-f^{(1)}$ solves
\begin{align*}
  \partial_t f_i(t)&=A_i^{(2)} \left(f_i^{(2)}(t)\right)^2 -A_i^{(1)}
\left(f_i^{(1)}(t)\right)^2+ B_i^\top
  f  +C^{(2)}_i(t)-C^{(1)}_i(t)\\
&=\left(A_i^{(2)} -A_i^{(1)} \right)\left(f_i^{(2)}(t)\right)^2
+A_i^{(1)} \left(f_i^{(2)}(t) +f_i^{(1)}(t)\right)f_i(t) + B_i^\top
  f (t)   +C^{(2)}_i(t)-C^{(1)}_i(t)\\
  &= \widetilde{B_i}(t)^\top f(t) + \widetilde{C_i}(t),
\end{align*}
where we write
\begin{align*}
  \widetilde{B_i}(t)&=B_i+ A_i^{(1)} \left(f_i^{(2)}(t)
+f_i^{(1)}(t)\right) e_i,\\ \widetilde{C_i}(t)&=\left(A_i^{(2)}
-A_i^{(1)}
\right)\left(f_i^{(2)}(t)\right)^2+C^{(2)}_i(t)-C^{(1)}_i(t).
\end{align*}
Note that $\widetilde{B}=(\widetilde{B}_{i,j})$ and $\widetilde{C}$
satisfy the assumptions of Corollary~\ref{corlincomp} in lieu of $B$
and $C$, and $f(0)\in\R^m_+$. Hence Corollary~\ref{corlincomp}
implies $f(t)\in\R^m_+$ for all $t\in [0,\tau)$, as desired. The
last statement of the lemma follows by the variation of constants
formula for $f^{(1)}(t)$.
\end{proof}

After these preliminary comparison results for the Riccati
equation~\eqref{ricccomp}, we now can state and prove an important
result for the system of Riccati equations~\eqref{eqriccfulladmiss}.
The following is an essential ingredient of the proof of Theorem
\ref{thmextanaXXX}. It is inspired by the line of arguments in
Glasserman and Kim~\cite{Glasserman}.
\begin{lemma}\label{lemstarshapedEXT}
Let $\Dcal_\R$ denote the maximal domain for the system of Riccati
equations~\eqref{eqriccfulladmiss}. Let $(\tau,u)\in\Dcal_\R$. Then
\begin{enumerate}
\item\label{lemstarshapedEXT1} $\Dcal_\R(\tau)$ is star-shaped around
zero.

\item\label{lemstarshapedEXT2} $\theta^\ast=\sup\{\theta\ge 0\mid
\theta u\in\Dcal_\R(\tau)\}$ satisfies either $\theta^\ast=\infty$
or $\lim_{\theta\uparrow\theta^\ast} \|\psi_I(t,\theta u)\|=\infty$.
In the latter case, there exists some $x^\ast\in\R^m_+\times \R^n$
such that $\lim_{\theta\uparrow\theta^\ast} \phi(\tau,\theta
u)+\psi(\tau,\theta u)^\top x^\ast=\infty$.
\end{enumerate}
\end{lemma}

\begin{proof}
We first assume that the matrices $\alpha_i$ are block-diagonal,
such that $\alpha_{i,iJ}=0$, for all $i=1,\dots,m$.

Fix $\theta\in (0,1]$. We claim that $\theta u\in\Dcal_\R(\tau)$. It
follows by inspection that $f^{(\theta)}(t)=\frac{\psi_I(t,\theta
u)}{\theta}$ solves \eqref{ricccomp} with
\[ A_i^{(\theta)}  =\frac{1}{2}\theta\alpha_{i,ii},\quad B =\Bcal_{II}^\top,
\quad C^{(\theta)}_i(t)=\beta_{i,J}^\top\psi_J(t,u)+
\frac{1}{2}\psi_J(t,u)^\top\theta\alpha_{i,JJ}\psi_J(t,u),\] and
$f(0)=u$. Lemma~\ref{lemcompriccX} thus implies that
$f^{(\theta)}(t)$ is nice behaved, as
\begin{equation}\label{lemstarshapedEXTeq1}
 \e^{\Bcal_{II}^\top t}\left(u+\int_0^t \e^{-\Bcal_{II}^\top
s}C^{(0)}(s)\,ds\right)\preceq f^{(\theta)}(t)\preceq \psi_I(t,u),
\end{equation}
for all $t\in [0,t_+(\theta u))\cap [0,\tau]$. By the maximality of
$\Dcal_\R$ we conclude that $\tau<t_+(\theta u)$, which implies
$\theta u\in\Dcal_\R(\tau)$, as desired. Hence $\Dcal_\R(\tau)$ is
star-shaped around zero, which is part~\ref{lemstarshapedEXT1}.

Next suppose that $\theta^\ast<\infty$. Since $\Dcal_\R(\tau)$ is
open, this implies $\theta^\ast u\notin\Dcal_\R(\tau)$ and thus
$t_+(\theta^\ast u)\le\tau$. From part~\ref{lemstarshapedEXT1} we
know that $(t,\theta u)\in\Dcal_\R$ for all $t<t_+(\theta^\ast u)$
and $0\le \theta\le \theta^\ast$. On the other hand, there exists a
sequence $t_n\uparrow t_+(\theta^\ast u)$ such that $
\|\psi_I(t_n,\theta^\ast u)\|>n$ for all $n\in\N$. By continuity of
$\psi$ on $\Dcal_\R$, we conclude that there exists some sequence
$\theta_n\uparrow\theta^\ast$ with $\|\psi_I(t_n,\theta_n
u)-\psi_I(t_n,\theta^\ast u)\|\le 1/n$ and hence
\begin{equation}\label{lemstarshapedEXTeq2}
 \lim_n\|\psi_I(t_n,\theta_n u)\|=\infty.
\end{equation}
Applying Lemma~\ref{lemcompriccX} as above, where initial time $t=0$
is shifted to $t_n$, yields
\[ g_n:=\e^{\Bcal_{II}^\top (\tau-t_n)}\left(f^{(\theta_n)}(t_n)
+\int_{t_n}^\tau \e^{\Bcal_{II}^\top
(t_n-s)}C^{(0)}(s)\,ds\right)\preceq f^{(\theta_n)}(\tau).\]
Corollary~\ref{corlincomp} implies that $\e^{\Bcal_{II}^\top
(\tau-t_n)}$ is $\succeq$-order preserving. That is,
$\e^{\Bcal_{II}^\top (\tau-t_n)}\R^m_+\subseteq\R^m_+$. Hence, in
view of \eqref{lemstarshapedEXTeq1} for $f^{(\theta_n)}(t_n)$,
\begin{align*}
g_n&\succeq \e^{\Bcal_{II}^\top (\tau-t_n)} \left(
\e^{\Bcal_{II}^\top t_n}\left(u+\int_0^{t_n} \e^{-\Bcal_{II}^\top
s}C^{(0)}(s)\,ds\right)+ \int_{t_n}^\tau \e^{\Bcal_{II}^\top
(t_n-s)}C^{(0)}(s)\,ds\right)\\ &=\e^{\Bcal_{II}^\top  \tau } \left(
u+\int_0^{\tau} \e^{-\Bcal_{II}^\top s}C^{(0)}(s)\,ds\right).
\end{align*}
On the other hand, elementary operator norm inequalities yield
\[ \|g_n\|\ge \e^{-\|\Bcal_{II}\|\tau} \|f^{(\theta_n)}(t_n)\| -
\e^{\|\Bcal_{II}\|\tau}\tau\sup_{s\in [0,\tau]}\|C^{(0)}(s)\|.\]
Together with \eqref{lemstarshapedEXTeq2}, this implies
$\|g_n\|\to\infty$. From Lemma~\ref{comppreceq} below we conclude
that $\lim_n f^{(\theta_n)}(\tau)^\top y^\ast =\infty$ for some
$y^\ast\in\R^m_+$. Moreover, in view of Lemma~\ref{lemcompriccX}, we
know that $f^{(\theta)}(\tau)^\top y^\ast$ is increasing $\theta$.
Therefore $\lim_{\theta\uparrow\theta^\ast}f^{(\theta)}(\tau)^\top
y^\ast=\infty$. Applying \eqref{lemstarshapedEXTeq1} and
Lemma~\ref{comppreceq} below again, this also implies that
$\lim_{\theta\uparrow\theta^\ast}\|f^{(\theta)}(\tau)\|=\infty$. It
remains to set $x^\ast=(y^\ast,0)$ and observe that $b_I\in\R^m_+$
and thus
\[ \phi(\tau,\theta u)=\int_0^\tau \left(\frac{1}{2}\psi_J(t,\theta
u)^\top a_{JJ}\,\psi_J(t,\theta u)+b_I^\top\psi_I(t,\theta
u)+b_J^\top\psi_J(t,\theta u)\right)dt \] is uniformly bounded from
below for all $\theta\in [0,\theta^\ast)$. Thus the lemma is proved
under the premise that the matrices $\alpha_i$ are block-diagonal
for all $i=1,\dots,m$.

The general case of admissible parameters $a,\alpha_i,b,\beta_i$ is
reduced to the preceding block-diagonal case by a linear
transformation along the lines of Lemma~\ref{lemcanaffinetrans}.
Indeed, define the invertible $d\times d$-matrix $\Lambda$
\[ \Lambda=\left(\begin{array}{cc} I_m & 0 \\
 {D} & I_n\end{array}\right)\]
where the $n\times m$-matrix ${D}=(\delta_1,\dots,\delta_m)$ has
$i$-th column vector
\[ \delta_{i}=\begin{cases} -\frac{\alpha_{i,iJ}}{\alpha_{i,ii}},&\text{if
$\alpha_{i,ii}>0$}\\ 0,&\text{else.}\end{cases}\] It is then not
hard to see that $\Lambda(\R^m_+\times\R^n)=\R^m_+\times\R^n$, and
\[\widetilde{\phi}(t,u)=\phi(t,\Lambda^\top u),\quad
\widetilde{\psi}(t,u)=\left(\Lambda^\top\right)^{-1}\psi(t,\Lambda^\top
u)\] satisfy the system of Riccati
equations~\eqref{eqriccfulladmiss} with $a,\alpha_i,b$, and
$\Bcal=(\beta_1,\dots,\beta_d)$ replaced by the admissible
parameters
\[ \widetilde{a}=\Lambda a\Lambda ^\top,\quad\widetilde{\alpha_i}=\Lambda\alpha_i\Lambda
^\top,\quad\widetilde{b}=\Lambda b,\quad\widetilde{\Bcal}=\Lambda
\Bcal\Lambda ^{-1}.\] Moreover, $\widetilde{\alpha_i}$ are
block-diagonal, for all $i=1,\dots,m$. By the first part of the
proof, the corresponding maximal domain
$\widetilde{\Dcal_\R}(\tau)$, and hence also
$\Dcal_\R(\tau)=\Lambda^\top\widetilde{\Dcal_\R}(\tau)$, is
star-shaped around zero. Moreover, if $\theta^\ast<\infty$, then
\[\lim_{\theta\uparrow\theta^\ast} \|\psi_I(\tau,\theta
u)\|=\lim_{\theta\uparrow\theta^\ast}
\left\|\widetilde{\psi}_I\left(\tau,\theta\left(\Lambda^\top\right)^{-1}
u\right)\right\|=\infty,\] and there exists some $ x^\ast
\in\R^m_+\times\R^n$ such that
\[ \lim_{\theta\uparrow\theta^\ast} \phi\left(\tau,\theta u\right)
+\psi\left(\tau,\theta u\right)^\top x^\ast
=\lim_{\theta\uparrow\theta^\ast}
\widetilde{\phi}\left(\tau,\theta\left(\Lambda^\top\right)^{-1}
u\right)
+\widetilde{\psi}\left(\tau,\theta\left(\Lambda^\top\right)^{-1}
u\right)^\top\Lambda x^\ast=\infty.\] Hence the lemma is proved.
\end{proof}

\begin{lemma}\label{comppreceq}
Let $c\in\R^m$, and $(c_n)$ and $(d_n)$ be sequences in $\R^m$ such
that
\[ c\preceq c_n\preceq d_n \]
for all $n\in\N$. Then the following are equivalent
\begin{enumerate}
  \item\label{comppreceq1} $\|c_n\|\to\infty$
  \item\label{comppreceq2} $  c_n^\top y^\ast\to\infty$ for some $y^\ast\in\R^m_+\setminus\{0\}$.
\end{enumerate}
In either case, $\| d_n\|\to\infty$ and $ d_n^\top y^\ast\to\infty$.
\end{lemma}

\begin{proof}
$\ref{comppreceq1}\Rightarrow\ref{comppreceq2}$: since
$\|c_n\|^2=\sum_{i=1}^m (  c_n^\top e_i)^2$ and $c_n^\top e_i\ge
c^\top e_i$, we conclude that $c_n^\top e_i\to\infty$ for some
$i=1,\dots,m$.

$\ref{comppreceq2}\Rightarrow\ref{comppreceq1}$: this follows from
$\|c_n^\top y^\ast\|\le  \|c_n\|\|y^\ast\|$.

The last statement now follows since $ d_n^\top y^\ast\ge  c_n^\top
y^\ast$.
\end{proof}

Finally, we sketch an alternative proof of Theorem
\ref{thmextanaXXX} \eqref{thmextanaXXX1} which avoids probabilistic
arguments.
\begin{remark}\label{alternative proof}
We may without loss of generality assume block-diagonal form of
$\alpha_i$, $i=1,\dots,d$ (cf.\ the final part of the proof of Lemma
\ref{lemstarshapedEXT}). Assume, by contradiction, that for some
$v\in\mathbb R^d$, $t_+(u+iv)<t_+(u)$. Then, as in the first proof,
we may deduce the existence of $t_n\uparrow t_+(u+iv)$ such that
\begin{equation}\label{extanaeq2x}
  \lim_{n}  (\Re \psi_i(t_n,u+iv))^+=\infty.
\end{equation}
 holds for some $i\in\{1,\dots,m\}$. Set $g(t,u+iv):=\Re(\psi_t,u+iv)$,
 $h:=\Im(\psi(t,u+iv)$. Then for $i=1,\dots,m$ the following
 differential inequality holds,
\begin{align}
\dot g_i(t,u+iv)&=\frac{1}{2}\alpha_{i,ii}(g_i^2-h_i^2) +g_J^\top
\alpha_{i,JJ}g_J-h_J^\top \alpha_{i,JJ}h_J+\beta_i^\top
g\\\nonumber&\leq\frac{1}{2}\alpha_{i,ii} g_i^2 +g_J^\top
\alpha_{i,JJ}g_J+\beta_i^\top g
\end{align}
and $g(t=0,u+iv))=\psi(t=0,u)=u$. Hence noting
$g_J(t,u+iv)=\psi_J(t,u)$ we obtain by Lemma \ref{lemcompriccX} for
all $t\in(0,t_+(u+iv))$
\[
\Re\psi(t,u+iv)=g(t,u+iv)\leq \psi(t,u).
\]
On the other hand, $\psi_I(t,u)\leq M$ for some positive constant
$M\in\mathbb R_+^m$, for all $t\in[0,t_+(u+iv)]$, hence
$\Re\psi_i(t,u+iv)\leq M_i$, which contradicts \eqref{extanaeq2x}.
\end{remark}
\end{appendix}

\bibliographystyle{siam}

\def\polhk#1{\setbox0=\hbox{#1}{\ooalign{\hidewidth
  \lower1.5ex\hbox{`}\hidewidth\crcr\unhbox0}}}

\end{document}